\documentclass[reqno,12pt,twoside,english]{amsart}
\usepackage[T1]{fontenc}
\usepackage[latin9]{inputenc}
\usepackage{amsthm}
\usepackage{amssymb}
\usepackage{amsmath,amsfonts,epsfig}

\usepackage{color}
\usepackage[final,linkcolor = blue,citecolor = blue,colorlinks=true]{hyperref}

\newtheorem*{thm*}{Theorem}

\newtheorem{thm}{Theorem}[section]
\newcommand{\bt}{\begin{thm}}
\newcommand{\et}{\end{thm}}

\newtheorem{cor}[thm]{Corollary}
\newcommand{\bc}{\begin{cor}}
\newcommand{\ec}{\end{cor}}

\newtheorem{lem}[thm]{Lemma}
\newcommand{\bl}{\begin{lem}}
\newcommand{\el}{\end{lem}}

\newtheorem{prop}[thm]{Proposition}
\newcommand{\bp}{\begin{prop}}
\newcommand{\ep}{\end{prop}}

\newtheorem{defn}[thm]{Definition}
\newcommand{\bd}{\begin{defn}}      
\newcommand{\ed}{\end{defn}}

\newtheorem{rmrk}[thm]{Remark}
\newcommand{\br}{\begin{rmrk}}
\newcommand{\er}{\end{rmrk}}

\newtheorem{quest}[thm]{Question}
\newcommand{\bq}{\begin{quest}}
\newcommand{\eq}{\end{quest}}

\newtheorem{example}[thm]{Example}

\newcommand\numberthis{\addtocounter{equation}{1}\tag{\theequation}}

\newcommand{\N}{\mathbb{N}}

\newcommand{\R}{\mathbb{R}}

\newdimen\vintkern\vintkern12pt
\def\vint{-\kern-\vintkern\int}

\newcommand{\B}{{\mathcal B}}
\newcommand{\Cp}{{\mathcal C}}
\newcommand{\Ch}{{\textbf {Ch}}}

\newcommand{\diam}{\operatorname{diam}}

\DeclareMathOperator{\loc}{loc}

\DeclareMathOperator{\modulus}{Mod}
\DeclareMathOperator{\capacity}{Cap}
\DeclareMathOperator{\supp}{supp}

\DeclareMathOperator{\RCD}{RCD}

\newcommand{\lip}{\operatorname{Lip}}

\newcommand*{\bR}{\ensuremath{\mathbb{R}}}

\def\Xint#1{\mathchoice
	{\XXint\displaystyle\textstyle{#1}}%
	{\XXint\textstyle\scriptstyle{#1}}%
	{\XXint\scriptstyle\scriptscriptstyle{#1}}%
	{\XXint\scriptscriptstyle\scriptscriptstyle{#1}}%
	\!\int}
\def\XXint#1#2#3{{\setbox0=\hbox{$#1{#2#3}{\int}$}
		\vcenter{\hbox{$#2#3$}}\kern-.5\wd0}}

\def\dashint{\Xint-}
\numberwithin{equation}{section}

\dedicatory{Dedicated to Professor~Pekka Koskela on the occasion of his $60^\text{th}$ birthday}

\begin{document}
\bibliographystyle{plain}

\title[Harmonic mappings between singular metric spaces]{Harmonic mappings between singular metric spaces}

\author[C.-Y. Guo]{Chang-Yu Guo}
\address[Chang-Yu Guo]{Research Center for Mathematics and Interdisciplinary Sciences, Shandong University 266237,  Qingdao, P. R. China}
\email{changyu.guo@sdu.edu.cn}

\subjclass[2010]{58E20,46E35,53C23,31C25}
\keywords{Dirichlet problem, Sobolev space, harmonic mapping, harmonic mapping flow, NPC space, ultra-completion}

\date{\today}

\thanks{C.-Y. Guo was supported by Swiss National Science Foundation Grant 175985 and the Qilu funding of Shandong University (No. 62550089963197).}

\begin{abstract}
In this paper, we survey the existence, uniqueness and interior regularity of solutions to the Dirichlet problem associated to various energy functionals in the setting of mappings between singular metric spaces. Based on known ideas and techniques, we separate the necessary analytical assumptions to axiomatizing the theory in the singular setting. More precisely,
\begin{itemize}
\item We extend the existence result of Guo and Wenger \cite{gw17} for solutions to the Dirichlet problem of Korevaar-Schoen energy functional to more general energy functionals in purely singular setting.

\item When $Y$ has non-positive curvature in the sense of Alexandrov (NPC), we show that the ideas of Jost \cite{j97} and Lin \cite{lin97} can be adapted to the purely singular setting to yield local H\"older continuity of solutions of the Dirichlet problem of Korevaar-Schoen and Kuwae-Shioya. 

\item We extend the Liouville theorem of Sturm \cite{s94} for harmonic functions to harmonic mappings between singular metric spaces.

\item We extend the theorem of Mayer \cite{m98} on the existence of the harmonic mapping flow and solve the corresponding initial boundary value problem.
\end{itemize}
Combing these known ideas, with the more or less standard techniques from analysis on metric spaces based on upper gradients, leads to new results when we consider harmonic mappings from $\RCD(K,N)$ spaces into NPC spaces. Similar results for the Dirichlet problem associated to the Kuwae-Shioya energy functional and the upper gradient functional are also derived.	

%One advantage of this type of axiomatization is  that it works for minimizers of other Dirichlet energy functional. In particular, as applications of the established theory, we deduce 

\end{abstract}

\maketitle
\tableofcontents

\section{Introduction}

Given a mapping $u\colon M\to N$ between two smooth Riemannian manifolds, there is a natural concept of energy associated to $u$. The minimizers, or more generally, the critical points of such energy functional, are called harmonic mappings. In the very beginning, the research on harmonic mappings comes together with the theory of minimal surfaces and it has attracted great attention after the work of Bochner~\cite{sy97}. However, the important existence, uniqueness and regularity theory were established relatively late - only after the work of Morrey~\cite{m48} on the Plateau problem in Riemannian manifold. The breakthrough in higher dimensional theory of harmonic mappings was made by Eells and Sampson~\cite{es64}, Hartman~\cite{h67} and idenpendently by Alber~\cite{alber64,alber67} and by Hamilton~\cite{h75} for manifolds with boundary via the heat equation method, where the target manifold $N$ was assumed to be non-positively curved. The regularity theory for general target Riemannian manifold has later been developed by Schoen and Uhlenbeck in a seminal paper~\cite{su82} and independently in~\cite{gg82,gg84} for manifolds with a single chart; see also~\cite{su83,hkw77,jm83}. In the remarkable work of Gromov and Schoen~\cite{gs92}, the authors proposed a variational approach for the theory of harmonic mappings to the setting of mappings into singular metric spaces, along with important applications to rigidity problems for certain discrete groups. 

Now, consider a mapping $u\colon X\to Y$, where $X=(X,d,\mu)$ is a metric measure space and $Y=(Y,d)$ a metric space. In a fundamental and important paper of Korevaar and Schoen~\cite{ks93}, an energy functional associated to $L^2(X,Y)$ mappings was introduced. More precisely, for each $\varepsilon>0$, one defines an approximating energy functional $E_\varepsilon(u)\colon C_0(X)\to \R$ on the space of continuous functions with compact support by
\begin{align*}
E_\varepsilon(u)(f)=\int_{X}f(x)\dashint_{B(x,\varepsilon)}\frac{d(u(x),u(y))^2}{\varepsilon^2}d\mu(y)d\mu(x).
\end{align*} 
In case $X$ is a compact $C^2$-smooth Riemannian manifold or a relatively compact domain in a $C^2$-smooth Riemannian manifold, it was shown that $E_\varepsilon(u)$ converges weakly, as a positive linear functional on $C_0(X)$, to some energy functional $E(u)$. Based on this energy functional, which we refer as Korevaar-Schoen energy functional, they have successfully extended the theory of harmonic mappings from $C^2$-smooth Riemannian manifolds into metric spaces with non-positive curvature in the sense of Alexandrov (NPC). Independently, Jost introduced in~\cite{j94} a slightly different energy functional and developed a theory of harmonic mappings associated to that energy functional through a sequential of deep works~\cite{j94,j95,j96,j97,j97book}. Moreover, the existence result of Jost~\cite{j94} works for mappings defined on more general metric spaces than domains in $C^2$-smooth Riemannian manifolds. From now on, if not specified, harmonic mappings refer to the energy minimizers of the Korevaar-Schoen energy functional.

Since then, there has been a considerable amount of growing interest in the theory of harmonic mappings between singular spaces. In particular, in the research monograph of Eells-Fuglede~\cite{ef01}, the authors extended the theory of harmonic mappings to the setting where $X$ is an admissible Riemannian polyhedron. Gregori~\cite{g98} further extended the existence and uniqueness theory of harmonic mappings to the setting where $X$ is a Lipschitz Riemannian manifold. Capogna and Lin~\cite{cl01} extended part of the  harmonic mapping theory to the setting of mappings from Euclidean spaces to the Heisenberg groups. In a series of deep works~\cite{st01,s02,s05}, Sturm developed a theory of harmonic mappings (associated to a slightly different functional) via a probabilistic theory and the theory of (generalized) Dirichlet forms.

\subsection{Existence and uniqueness}
There are two general approaches for the existence of harmonic mappings: the first one relies on the uniform convexity of the distance function in NPC spaces (see e.g.~\cite{ks93,g98,ef01,f05}), while the second one relies  on the theory of (metric space valued) Sobolev mappings, in particular, the theory of trace, lower semicontinuity of the enery functional with respect to the $L^2$-convergence and the (various versions of) Rellich compactness theorem (see e.g.~\cite{lw16,gw17}). In the first approach, a crucial fact one needed is the so-called subpartition lemma, which essentially says that the integral averages one uses to approximate the Sobolev energy satisfy certain monotonicity with respect to the size of the ball on which the average is taken. The advantage of this approach is that we can solve the Dirichlet problem of Korevaar and Schoen for fairly general open subset of the metric measure space $X$; see for instance~\cite{f05}. In the second approach, the domain $\Omega\subset X$ has to be sufficiently nice so that both the theory of trace and certain version of Rellich compactness theorem holds. The advantage of the second approach is that the target metric space $Y$ does not need to NPC, and indeed, can be fairly general, which includes in particular all proper metric spaces, all dual Banach spaces and NPC spaces; see~\cite{gw17}. Of course, one cannot expect uniqueness in such a great generality.

Our first main result of this paper concerns the existence of harmonic mappings into singular metric spaces. Note that both our source domain and the target metric space are fairly general.

\bt\label{thm:main thm existence and uniqueness}
Let $(X,d,\mu)$ be a compact metric space and $Y$ a metric space that is 1-complemented in some ultra-completion of $Y$. Suppose $X$ admits an energy functional $\mathcal{E}$ with property $\B$. Fix a domain $\Omega\subset X$. 
%that supports a weak $(1,2)$-Poincar\'e inequality~\eqref{eq:Poin ineq}. 
Then, for each $\phi\in S^{1,2}(X,Y)$, there exists a mapping $u\in S^{1,2}_{\phi}(\Omega,Y)$ such that 
\begin{equation*}
  \mathcal{E}(u)=\inf_{v\in S^{1,2}_{\phi}(\Omega,Y)} \mathcal{E}(v).
\end{equation*}	
\et

As observed in~\cite[Proposition 2.1]{gw17}, many nice metric spaces are 1-complemented in some ultra-completion of themselves, in particular, proper metric spaces, dual Banach spaces, $L^1$-spaces and NPC spaces. 
The definition of metric spaces that admits an energy functional $\mathcal{E}$ with property $\B$  shall be given in Section~\ref{subsec:metric spaces with property B} below. It is stated in a very abstract sense. 

%However, one can verify that  the Korevaar-Schoen energy functional $E$, the Kuwae-Shioya energy functional $E^b$ and the upper gradient energy functional $E^g$ satisfy the assumptions of property $\mathcal{B}$. Consequently, as applications of Theorem \ref{thm:main thm existence and uniqueness}, we obtain the existence for the Dirichlet problem associated to these energy functionals.

As  the Korevaar-Schoen energy functional $E$ has property $\mathcal{B}$, we obtain, as a direct corollary of Theorem \ref{thm:main thm existence and uniqueness}, the existence of Dirichlet problem of Korevaar-Schoen.
\bc\label{coro:KS1}
Suppose $(X,d,\mu)$ is compact $\RCD(K,N)$ and $Y$ is a metric space that is 1-complemented in some ultra-completion of $Y$. Then for each $\phi\in KS^{1,2}(X,Y)$, there exists a mapping $u\in KS^{1,2}_{\phi}(\Omega,Y)$ such that 
\begin{equation*}
	E(u)=\inf_{v\in KS^{1,2}_{\phi}(\Omega,Y)} E(v).
\end{equation*}	
\ec

%Many nice spaces have strong property $\B$, in particular, Lipschitz manifolds considered in~\cite{g98}, admissible Riemannian polyhedrons considered in~\cite{c95,ef01,dm10}, metric spaces with the strong measure contraction property (SMCP) considered in~\cite{s98} and $\RCD(K,N)$-spaces considered in \cite{ags14b,ags15,gt20}. 

Kuwae and Shioya~\cite{ks03} constructed an energy functional which slightly differs from the Korevaar-Schoen energy functional. They introduced a new concept, called \emph{strong measure contraction property of Bishop-Gromov type (SMCPBG)} (given in Section~\ref{subsec:SMCP} below), and developed a theory of Sobolev space of mappings $u\colon X\to Y$ when the source metric measure space space $X$ satisfies the SMCPBG. Riemannian manifolds and Alexandrov spaces with curvature bounded from below are typical examples of metric measure spaces that posse the SMCPBG. Since the Kuwae-Shioya energy functional $E^b$ has property $\B$, we have the following corollary. 

\bc\label{coro:KS2}
Suppose $(X,d,\mu)$ is a compact metric space satisfying the SMCPBG and $Y$ is a metric space that is 1-complemented in some ultra-completion of $Y$. Then for each $\phi\in W^{1,2}(X,Y)$, there exists a mapping $u\in W^{1,2}_{\phi}(\Omega,Y)$ such that 
\begin{equation*}
	E^b(u)=\inf_{v\in W^{1,2}_{\phi}(\Omega,Y)} E^b(v).
\end{equation*}	
\ec

Parallel to the theory of harmonic mappings, great effort has been made to extend the existence, uniqueness and regularity theory of harmonic functions, that is harmonic mappings into the real line $Y=\R$, to the abstract metric measure space setting; see for instance~\cite{c99,s01,ks01,krs03,hkm06,j14,bb11} and the references therein. Unlike the case of mappings, one usually uses the $L^2$-norm of the \emph{upper gradients} as the energy functional. Since the notion of upper gradients works also for mappings, it is natural to consider the Dirichlet problem associated to the energy functional of upper gradients. As a by-product of Theorem~\ref{thm:main thm existence and uniqueness}, we obtain the following existence result for the Dirichlet problem based on upper gradients.

\bc\label{coro:UG}
Suppose $(X,d,\mu)$ is a compact PI space and $Y$ is a metric space that is 1-complemented in some ultra-completion of $Y$. Then for each $\phi\in N^{1,2}(X,Y)$, there exists a mapping $u\in N^{1,2}_{\phi}(\Omega,Y)$ such that 
\begin{equation*}
	E^g(u)=\inf_{v\in N^{1,2}_{\phi}(\Omega,Y)} E^g(v).
\end{equation*}	
\ec

Recall that we say a metric measure space $X=(X,d,\mu)$ a PI-space, if the measure $\mu$ is doubling on $X$, i.e., there exists a constant $c_d>0$ such that
$$\mu(B(x,2r))\leq c_d\mu(B(x,r))$$
for all open balls $B(x,r)\subset X$ with $\diam B\leq \diam X$ and it supports a weak (1,2)-Poincar\'e inequality~\eqref{eq:Poin ineq}. Corollary~\ref{coro:UG} can be regarded as a natural extension of~\cite[Theorem 5.6]{s01} to the setting of metric space valued mappings. Note however that the proof of Theorem 5.6 in~\cite{s01} cannot work in this case since $N^{1,2}(X,Y)$ is not a linear (Banach) space for general metric target.
%Thus Theorem~\ref{thm:main thm existence and uniqueness} provides a unified treatment for the existence and uniqueness theory of harmonic mappings in the singular setting.

The idea for the proof of Theorem~\ref{thm:main thm existence and uniqueness} relies on that used in~\cite[Proof of Theorem 1.4]{gw17} but some more care needs to be paid in order to deal with the singular/irregular source domain. The main difference with the situation there is that we do not have a theory of trace for metric space valued Sobolev mapping as that of~\cite{ks93}. Instead, we (essentially) use the fact that admissible mappings are defined on a (relatively compact) neighborhood of the domain $\Omega$ and coincide with a given boundary value. In the first step, we apply the generalized compactness result (see~\cite[Theorem 3.1]{gw17}) to obtain an energy minimizing harmonic mapping into some ultra-completion $Y_\omega$ of the target metric space $Y$. Composing this map with the 1-Lipschitz projection map from $Y_\omega$ to $Y$, we then obtain an energy minimizing harmonic mapping from $X$ to $Y$. The only remaining issue is to show that this map has the correct boundary value, which will follow by our construction of the map into the ultra-completion.
%This has the advantage that we can solve the Dirichlet problem for fairly general domains. 

Regarding the uniqueness, we have the following result for the Korevaar-Schoen energy functional. 
\bt\label{thm:uniqueness}
	Suppose $(X,d,\mu)$ is $\RCD(K,N)$ and $Y$ is NPC. Fix a relatively compact domain $\Omega\subset X$ that supports a $(1,2)$-Poincar\'e inequality~\eqref{eq:omega weak poincare} and that has the property $\capacity_2(X\backslash \Omega)>0$. Then, for each $\phi\in KS^{1,2}(X,Y)$, there exists a unique mapping $u\in KS^{1,2}_{\phi}(\Omega,Y)$ such that 
	\begin{equation*}
		E(u)=\inf_{v\in KS^{1,2}_{\phi}(\Omega,Y)} E(v).
	\end{equation*}	
\et

The general idea of the proof of Theorem~\ref{thm:uniqueness} is quite similar to~\cite[Proof of Theorem 2.2]{ks93}, but relies on the theory of metric-valued Sobolev spaces developed in~\cite{hkst12}. In particular, we provide a new and concrete proof of~\cite[Theorem 1 (a)]{f05} in a greater generality. Moreover, if $X$ is instead assumed to satisfy the SMCPBG, then the conclusion of Theorem \ref{thm:uniqueness} holds with the Kuwae-Shioya energy functional $E^b$; see Remark \ref{rmk:on existence and uniqueness}.

\subsection{Interior regularity}
In their fundamental work~\cite{ks93}, Korevaar and Schoen has shown that harmonic mappings from a $C^2$-smooth Riemannian manifold to an NPC space is locally Lipschitz continuous, which plays an important role in establishing rigidity theorems of geometric group theory (see e.g.~\cite{gs92,dm08}). Since then, there has been a lot of effort in establishing interior regularity of harmonic mappings in the singular space setting; see for instance~\cite{c95,s96,j95,j97,f03,f05,f08,iw08,dm08,zz13,h16,hz16,zz16,lw16} and the references therein. It should be noticed that one cannot expect local Lipschitz continuity holds in general singular metric spaces. Indeed, Chen~\cite{c95} constructed a harmonic function on a two-dimensional metric cone $X$ such that $u$ is not Lipschitz continuous if $X$ has no lower curvature bound. Nevertheless, harmonic functions constructed there do satisfy the local H\"older continuity, which is valid for all harmonic mappings from admissible Riemannian polyhedrons to NPC spaces (see e.g.~\cite{c95,ef01}). 

In~\cite{lin97}, Lin proposed a very elegant approach to obtain H\"older continuity of harmonic mappings between singular spaces. In particular, Lin's method implies that  harmonic mappings from Alexandrov spaces with curvature bounded from below to locally doubling NPC spaces are local H\"older continuous, provided that the composition of the distance function with the harmonic mapping is subharmonic. As observed in~\cite{zz13}, the latter requirement holds for harmonic mappings from Alexandrov spaces with curvature bounded from below to NPC spaces via an argument due to Jost~\cite{j97}. Note that in the work of Jost~\cite{j97}, a fairly general interior H\"older regularity result for harmonic mappings (associated to the Jost's energy functional) was established. On the other hand, as pointed out in~\cite{kst04}, the approach of Jost~\cite{j97} relies on a theory of generalized Dirichlet forms for metric space valued mappings, which is in general hard to verify since the classical method of Beurling and Deny (see e.g.~\cite{fot94}) fails in constructing energy measures for metric-space valued mappings. However, the essential ingredients in Jost's regularity result are volume doubling property and (2,2)-Poincar\'e inequality for the intrinsic metric space (induced by the Dirichlet form) and hence it can be extended to rather general setting.
%For this reason, it is not even clear whether one can apply the result of Jost~\cite{j97} to conclude the interior H\"older regularity of Korevaar-Schoen harmonic mappings in the setting of Lipschitz Riemannian manifolds as considered in~\cite{g98}.
%Korevaar-Schoen harmonic mappings are locally H\"older continuous in the setting of~\cite{g98}.

Our second main purpose of this paper is to establish interior H\"older regularity of Korevaar-Schoen and Kuwae-Shioya harmonic mappings in a large class of singular metric spaces. Our second main result reads as follows.

\bt\label{thm:main thm regularity}
1). Let $X$ be a complete metric measure space and $Y$ a locally doubling NPC space. If the approximating energy density $e_\varepsilon^{c,Y}$ between $X$ and $Y$ has property $\Cp$, then each solution $u$ of the Dirichlet problem (if exists) is locally H\"older continuous.

2). Let $X$ be a complete metric measure space and $Y$ an NPC space. If the approximating energy density $e_\varepsilon^{c,Y}$ between $X$ and $Y$ has strong property $\Cp$, then each solution $u$ of the Dirichlet problem (if exists) is locally H\"older continuous.
\et 

The definition of approximating energy density with (strong) property $\Cp$ shall be given in Section~\ref{subsec:metric spaces C} below. For a given NPC space $Y$, when $X$ is either the Lipschitz manifold considered in~\cite{g98}, or the admissible Riemannian polyhedron considered in~\cite{c95,ef01,dm10}, or metric spaces with the strong measure contraction property (SMCP) considered in~\cite{s98} or $\RCD(K,N)$-space considered in \cite{ags14b,ags15,gt20}, the Korevaar-Schoen approximating energy density between $X$ and $Y$ will have strong property $\Cp$. Similarly, when $X$ satisfies the SMCPBG, the Kuwae-Shioya approximating energy density will have strong property $\Cp$; see Examples \ref{ex:metric space C} and \ref{exam:Kuwae-Shioya strong property c} below. Consequently, we obtain the following corollary of Theorem \ref{thm:main thm regularity}.
%It includes many nice spaces, in particular, Lipschitz manifolds considered in~\cite{g98}, admissible Riemannian polyhedrons considered in~\cite{c95,ef01,dm10}, metric spaces with the strong measure contraction property (SMCP) considered in~\cite{s98} and $\RCD(K,N)$-spaces considered in \cite{ags14b,ags15,gt20}. 

\bc\label{coro:interior regularity}
Let $Y$ be an NPC space.  Then we have
\begin{itemize}
\item when $X$ is $\RCD(K,N)$, each solution $u$ of the Dirichlet problem of Korevaar-Schoen (if exists) is locally H\"older continuous;
\item when $X$ satisfies the SMCPBG, each solution $u$ of the Dirichlet problem of Kuwae-Shioya (if exists) is locally H\"older continuous.
\end{itemize}

\ec

The proof of Theorem~\ref{thm:main thm regularity} 1) is a combination of the approach of Jost~\cite{j97} and Lin~\cite{lin97}. More precisely, we first follow the idea of Jost to show that for each point $y_0\in Y$, the function $f_{y_0}:=d^2(u(\cdot),y_0)$ is (weakly) subharmonic in the sense of~\cite{bm95} (in terms of Dirichlet forms), and then adapt the argument of Lin~\cite{lin97} to prove the local H\"older regularity. We would like to point out that the argument of Lin~\cite{lin97} is elegant and beautiful, but it requires the target space to be (locally) doubling. Thus it cannot be applied in general NPC targets.

The proof of Theorem~\ref{thm:main thm regularity} 2) follows closely the approach of Jost~\cite{j97}. The essential difference with~\cite{j97} is that we work directly on the Dirichlet forms for functions, instead of the generalized Dirichlet forms for mappings as in~\cite{j97}, whereas the well-known regularity theory of sub/super solutions associated to Dirichlet forms developed in~\cite{bm95} can be applied. 

We did not address the interior regularity for the case when $Y$ is a CAT(1) space in this paper. It is worth pointing out that in this more general case similar interior regularity results have been obtained when $X$ is an admissible Riemannian polyhedron in~\cite{f03,f08} or when $X$ is an Alexandrov space with curvature bounded from below in~\cite{hz16,zzz19}.

Note that Theorem \ref{thm:main thm regularity} implies that harmonic mappings (associated to the Korevaar-Schoen energy functional) from $\RCD(K,N)$ space to NPC space are locally H\"older continuous. In a recent remarkable work of Zhang and Zhu~\cite{zz13} (see also~\cite{zz16}), harmonic mappings were shown to be locally Lipschitz continuous, when the source metric measure space is an Alexandrov space with curvature bounded from below and the target metric space is NPC. Since Alexandrov spaces with curvature bounded from below are special cases of $\RCD(K,N)$-spaces, this and the corresponding result for harmonic functions~\cite{j14} seem to suggest that Korevaar-Schoen harmonic mappings may be locally Lipschitz continuous, if the metric measure space $X$ is $\RCD(K,N)$; see~\cite[Problem 4.6]{zz16} for more detailed discussions on this interesting open problem.
%On the other hand, the notion of a lower bound on the Ricci curvature for general metric measure spaces has been introduced by Lott-Villani-Sturm~\cite{lv09,s06} and even a Riemannian Ricci lower bound, the so-called RCD$^*(K,N)$ spaces, was introduced in~\cite{ags14b,eks15,ags16}. Thus it would be very interesting to know whether one can further extend the work of Zhang and Zhu to the RCD$^*(K,N)$ spaces; see~\cite[Problem 4.6]{zz16} for more detailed discussions.

\textbf{Open problem.} Can we establish the interior H\"oler continuity of solutions for the Dirichlet problem associated to the upper gradient energy functional?

The main difficulty to this open problem is that it is very difficult to compute the upper gradient of a composed function and thus we cannot use the approach for Theorem \ref{thm:main thm regularity} to infer the composition of the distance function with $u$ will be any kind of ``subharmonic" function.

\subsection{Liouville theorem}

Liouville type theorems for harmonic mappings between complete smooth Riemannian manifolds have been investigated by many authors including geometers and probabilists. In particular, Eells-Sampson~\cite{es64} proved that any bounded harmonic mapping from a compact Riemannian manifold with positive Ricci curvature into a complete manifold with non-positive curvature is constant. Schoen and Yau~\cite{sy76} proved that any harmonic mapping with finite energy from a complete smooth Riemannian manifolds with non-negative Ricci curvature into a complete manifold with non-positive curvature is constant. Cheng~\cite{c80} showed any harmonic mapping with sublinear growth from a complete Riemannian manifold with non-negative Ricci curvature into an Hadamard manifold is constant. Hildebrandt-Jost-Widman~\cite{hjw80} proved a Liouville type theorem for harmonic mappings into regular geodesic balls in a complete smooth Riemannian manifold. For a detailed description of other types of Liouville type theorems for harmonic mappings; see~\cite{ks08}.

For the statement of our Liouville theorem for harmonic mappings, we set for $u\colon X\to Y$, 
$$v_p(r,x,x_0)=\int_{B(x,r)}d(u(x),u(x_0))^pd\mu(x)$$

\bt\label{thm:Liouville}
Let $X$ be a complete metric measure space and $Y$ an NPC space. Suppose the approximating energy density $e_\varepsilon^{c,Y}$ between $X$ and $Y$ has property $\Cp$. If $u\in KS^{1,2}_{\loc}(X,Y)$ is a harmonic mapping such that for some $x,x_0\in X$ and $p>1$, 
\begin{equation}\label{eq:divergence condition for Liouville}
	\int_1^\infty \frac{r}{v_p(r,x,x_0)}dr=\infty,
\end{equation}
then $u$ is constant.
\et

The proof of Theorem~\ref{thm:Liouville} relies on the Liouville type theorem for weakly subharmonic functions. Originally, Yau~\cite[Theorem 1]{y76} has shown that there is no non-constant smooth non-negative $L^p$-integrable, $p>1$, subharmonic functions on a complete smooth Riemannian manifold. Sturm~\cite[Theorem 1]{s94} extended this result to the setting of Dirichlet forms under the sharper condition~\eqref{eq:divergence condition for Liouville}. In our setting, Theorem~\ref{thm:Liouville} is a direct application of this Liouville theorem for weakly subharmonic functions. More precisely, we shall show that $v=d(u(\cdot),u(x_0))$ is a weakly subharmonic function (in the sense of Sturm~\cite{s94}) on $X$ and so by~\cite[Theorem 1]{s94}, $v$ is constant. Consequently, $u=u(x_0)$ is constant. In the proof of $v$ being weakly subharmonic, we shall combine the idea of Jost \cite{j97} together with the theory of fine topology and potential theory that developed recently in the setting of metric spaces \cite{bb11,bbl18}.

\subsection{Harmonic mapping flow}
In~\cite{m98}, Mayer developed a general theory of gradient flows on NPC spaces with successful applications to harmonic mapping. The basic setting is an NPC space $L=(L,D)$ together with an energy functional $G\colon L\to \R\cup\{\infty\}$. The gradient flow equation
\begin{align*}
\frac{du(t)}{dt}=-\nabla G(u(t))
\end{align*}
has been interpreted as a variational formulation:
\begin{align*}
u(t+h) \quad  \text{minimizes}\quad  u\mapsto G(u)+\frac{1}{2h}D^2(u,u(t)).
\end{align*}
The existence and uniqueness of $u(t+h)$ were obtained under mild structural assumptions on $G$; see~\cite[Theorem 1.13]{g98} and \cite{j98}.

In our setting, it is not difficult to show that $L=L^2(X,Y)$ is NPC when $X$ is complete and $Y$ is NPC (see Lemma~\ref{lemma:L2 is NPC}). On the other hand, when $X$ has strong property $\B$, $u\mapsto E(u)$ is a lower semicontinuous convex functional. Thus as an immediate corollary of~\cite[Theorem 1.13]{m98}, we obtain the following result.
%Applying the method of Mayer~\cite{m98}, we deduce the following result.

\bt[Existence and boundedness of gradient flow of the Dirichlet Energy]\label{thm:gradient flow}
Assume $X$ is $\RCD(K,N)$ and $Y$ is NPC. For any starting point $u_0\in KS^{1,2}(X,Y)$ the gradient flow of the Korevaar-Schoen energy exists and $u_t\in KS^{1,2}(X,Y)$ for each $t\geq 0$. Moreover, if $X$ has finite $\mu$-measure, then the gradient flow stays bounded for all times.
\et 

As an immediate corollary of Theorem~\ref{thm:gradient flow} and the Rellich compactness theorem for Sobolev mappings, we obtain the following corollary.
\bc\label{coro:gradient flow}
In the setting of Theorem~\ref{thm:gradient flow}, if $X$ is compact and $Y$ is proper, then the flow $u_t$ converges to a constant mapping as $t\to\infty$.
\ec 

As an easy consequence of our proof of Theorem~\ref{thm:main thm existence and uniqueness}, we can show that the Sobolev space $\mathcal{KS}^{1,2}_\phi(\Omega,Y)$ is NPC and consequently solve the initial boundary value problem via the exactly same method as in~\cite[Theorem 3.4]{m98}.

\bt[Solvability of the initial boundary value problem]\label{thm:initial boundary value problem}
Assume $X$ is $\RCD(K,N)$ and $Y$ is NPC. For any given map $\phi\in KS^{1,2}(\Omega,Y)$, the following problem admits a solution in the sense of~\cite[Theorem 1.13]{m98}:
\begin{equation*}
\begin{cases} u(t) \text{ solves the harmonic mapping flow for } t\geq 0, \\
u(0)=\phi,\\
\text{a representative of } u(t) \text{ equal }\phi \text{ q.e. in }X\backslash \Omega.\\
\end{cases}
\end{equation*}
Moreover, if $\Omega$ is relatively compact in $X$, then $\hat{u}=\lim_{t\to\infty}u(t)$ exists and is the unique harmonic mapping solving the Dirichlet problem of Korevaar and Schoen with boundary data $\phi$.
\et

\textbf{Sturcture of the paper.}
This paper is structured as follows. In section~\ref{sec:preliminaries}, we recall the necessary definitions needed for this paper, mainly, the definition of Sobolev spaces of metric-valued mappings. In section~\ref{sec:existence and uniqueness for KS}, we prove the existence and uniqueness theorem. In section~\ref{sec:holder continuity of KS}, we show the interior H\"older continuity of harmonic mappings. In section~\ref{sec:Liouville}, we prove the Liouville theorem for harmonic mappings. We study the associated harmonic mapping flow in Section~\ref{sec:harmonic mapping flow}. In the final section, Section~\ref{sec:Dirichlet problem for kuwae-Shioya}, we study the Dirichlet problem associated to other energy functionals, in particular, the energy functional of Kuwae and Shioya and the energy functional based on upper gradients.

\section{Preliminaries}\label{sec:preliminaries}

The metric space $Y$ in this paper is always assumed to be separable so that we may embeds $Y$ isometrically into the Banach space $l^\infty(Y)$.

\subsection{Metric spaces of non-positive curvature}

We shall need the following concept, which was introduced by Alexandrov \cite{Ale57}; see also \cite{bn93}.
\bd[NPC spaces]
A complete metric space $(X,d)$
(possibly infinite dimensional) is said to be non-positively curved
(NPC) if the following two conditions are satified:

\begin{itemize}
	\item $(X,d)$ is a length space, that is, for any two points $P,Q$ in
	$X$, the distance $d(P,Q)$ is realized as the length of a rectifiable
	curve connecting $P$ to $Q$. (We call such distance-realizing curves
	geodesics.)
	\item For any three points $P,Q,R$ in $X$ and choices of geodesics $\gamma_{PQ}$
	(of length $r$), $\gamma_{QR}$ (of length $p$), and $\gamma_{RP}$ (of
	length $q$) connecting the respective points, the following comparison
	property is to hold: For any $0<\lambda<1$, write $Q_{\lambda}$ for the
	point on $\gamma_{QR}$ which is a fraction $\lambda$ of the distance from
	$Q$ to $R$. That is, 
	\[
	d(Q_{\lambda},Q)=\lambda p,\quad d(Q_{\lambda},R)=(1-\lambda)p.
	\]
	On the (possibly degenerate) Euclidean triangle of side lengths $p,q,r$
	and opposite vertices $\bar{P}$, $\bar{Q},\bar{R}$, there is a corresponding
	point 
	\[
	\bar{Q}_{\lambda}=\bar{Q}+\lambda(\bar{R}-\bar{Q}).
	\]
	The NPC hypothesis is that the metric distance $d(P,Q_{\lambda})$ (from
	$Q_{\lambda}$ to the opposite vertex $P$) is bounded above by the Euclidean
	distance $|\bar{P}-\bar{Q}_{\lambda}|$. This inequality can be written
	precisely as 
	\[
	d^{2}(P,Q_{\lambda})\le(1-\lambda)d^{2}(P,Q)+\lambda d^{2}(P,R)-\lambda(1-\lambda)d^{2}(Q,R).
	\]
\end{itemize}
\ed 

In an NPC space $Y$, geodesics connecting each pair of points are unique and so one can define the $t$-fraction mapping $u_t$ of two mapping $u_0,u_1\colon X\to Y$ as $u_t=``(1-t)u_0+tu_1"$, that is, for each $x$, $u_t(x)$ is the unique point $P$ on the geodesic connecting $u_0(x)$ and $u_1(x)$ such that $d(P,u_0(x))=td(u_0(x),u_1(x))$ and $d(P,u_1(x))=(1-t)d(u_0(x),u_1(x))$. We refer the interested readers to~\cite[Section 2.1]{ks93} or~\cite{j97book} for more discussions on NPC spaces.

\subsection{Strong measure contraction properties}\label{subsec:SMCP}
The following notion of measure contraction property was introduced by Sturm~\cite{s98}.
\bd[Weak measure contraction property]\label{def:weak measure contraction property}
We say that a metric measure space $(X,d,\mu)$ satisfies the weak measure contraction property (WMCP) if there exist numbers $R>0$, $\theta<\infty$ and $\lambda<\infty$ and $\mu^2$-measurable maps $\Phi_t\colon X\times X\to X$ with the following properties:

(1) For $\mu$-a.e. $x_1,x_2\in X$ with $d(x_1,x_2)<R$ and all $s,t\in [0,1]$, 
\begin{equation}\label{eq:MCP 1}
	\Phi_0(x,y)=x,\Phi_t(x,y)=\Phi_{1-t}(y,x),\Phi_{s}(x,\Phi_t(x,y))=\Phi_{st}(x,y),
\end{equation}
and
\begin{equation}\label{eq:MCP 2}
	d(\Phi_s(x,y),\Phi_t(x,y))\leq \lambda |s-t|d(x,y)
\end{equation}

(2) For all $\varepsilon<R$, $\mu$-a.e. $x\in X$, all $\mu$-measurable set $A\subset B(x,r)$ and all $t\in [0,1]$, 
\begin{equation}\label{eq:MCP 3}
	\frac{\mu_\varepsilon(A)}{\sqrt{\mu(B(x,r))}}\leq \theta \frac{\mu_{t\varepsilon}(\Phi_t(x,A))}{\sqrt{\mu(B(x,t\varepsilon))}},
\end{equation}
where $d\mu_\varepsilon(x)=\frac{d\mu(x)}{\sqrt{\mu(B(x,\varepsilon))}}$ and $d\mu_\varepsilon(y)=\frac{d\mu(y)}{\sqrt{\mu(B(y,\varepsilon))}}$.
\ed

\bd[SMCP]\label{def:SMCP}
We say that a metric measure space $(X,d,\mu)$ posses the strong measure contraction property (SMCP) if it satisfies the WMCP and the constants $\Theta$ and $\theta$ appearing in Definition~\ref{def:weak measure contraction property} can be chosen to be arbitrarily close to 1.
\ed 

Many interesting metric spaces satisfies SMCP. In particular, $C^2$-smooth Riemannian $n$-manifolds; see~\cite[Section 4]{s98} for more examples.

\bd[WMCPBG]\label{def:WMCPBG}
We say that a metric measure space $(X,d,\mu)$ satisfies the weak measure contraction property of Bishop-Gromov type (WMCPBG) if there exist $\nu>0$, and increasing continuous function $b\colon [0,\infty)\to [0,\infty)$ with $b(0)=0$ such that there are positive finite constants $R$, $\Theta\geq 1$ and $\theta\geq 1$, and $\mu\otimes \mu$-measurable mappings $\Phi_t=\Phi_t^{Z}\colon X\times X\to X$ for all $t\in [0,1]$, with the following
\begin{itemize}
	\item[i).] For $\mu$-a.e. $x,y\in X$ with $d(x,y)<R$ and all $s,t\in [0,1]$,
	\begin{align*}
		\Phi_0(x,y)=x, &\Phi_t(x,y)=\Phi_{1-t}(y,x), \Phi_s(x,\Phi_t(x,y))=\Phi_{st}(x,y)\\
		&d(\Phi_s(x,y),\Phi_t(x,y))\leq \theta|t-s|d(x,y).
	\end{align*} 
	
	\item[ii).] For all positive $r<R$, $\mu$-a.e. $x\in X$, all $\mu$-measurable set $A\subset B(x,r)$ and all $t\in [0,1]$,
	\begin{align*}
		\frac{\mu(A)}{b(r)}\leq \Theta\frac{\mu(\Phi_t(x,A))}{b(rt)}.
	\end{align*}

	\item[iii).] For all positive $r<R$, $\mu$-a.e. $x\in X$,
	\begin{align*}
		\mu(B(x,r))\leq \theta b(r).
	\end{align*}
	
	\item[iv).] For $0<r_1\leq r_2\leq R\theta^2$,
	\begin{align*}
		\frac{b(r_2)}{b(r_1)}\leq \Theta\Big(\frac{r_2}{r_1}\Big)^\nu.
	\end{align*}	
\end{itemize}
\ed 

\bd[SMCPBG]\label{def:SMCPBG}
We say that a metric measure space $(X,d,\mu)$ posses the strong measure contraction property of Bishop-Gromov type (SMCPBG) if it satisfies the WMCPBG and the constants $\Theta$ and $\theta$ appearing in Definition~\ref{def:WMCPBG} can be chosen to be arbitrarily close to 1.
\ed

\begin{example}\label{example:SMCPBG}
	The following spaces satisfy the SMCPBG property.	
	\begin{itemize}
		\item[1).] $C^2$-smooth Riemannian manifolds $(M,d_g,\mu_g)$ satisfies SMCPBG with $b(r)=\omega_{n-1}r^n/n$.
		
		\item[2).] Let $X$ be an $n$-dimensional Alexandrov space with curvature bounded from below by some $\kappa\in \R$. Then $X$ has the SMCPBG with $\mu=\mathcal{H}^n$ and $b(r)=\omega_{n-1}r^n/n$. 
	\end{itemize}
\end{example}

More interesting examples can be found in~\cite[Section 4]{s98}.

\subsection{Ultra-completions of metric spaces}
We briefly recall the relevant definitions concerning ultra-completions and ultra-limits of metric spaces. Details can be found for instance in \cite{bh99}.

A non-principal ultrafilter on $\N$ is a finitely additive probability measure $\omega$ on $\N$ such that every subset of $\N$ is measurable and such that $\omega(A)$ equals 0 or 1 for all $A\subset \N$ and $\omega(A)=0$ whenever $A$ is finite. Given a compact Hausdorff topological space $(Z,\tau)$ and a sequence $\{z_m\}\subset Z$ there exists a unique point $z_\infty\in Z$ such that $\omega(\{m\in \N:z_m\in U\})=1$ for every $U\ni \tau$ containing $z_\infty$. We denote the point $z_\infty$ by $\lim_\omega z_m$.

Let $Y=(Y,d)$ be a metric space and $\omega$ a non-principal ultrafilter on $\N$. A sequence $\{y_m\}\subset Y$ is bounded if $\sup_m d(y_1,y_m)<\infty$. Define an equivalence relation $\sim$ on bounded sequences in $Y$ by considering $\{y_m\}$ and $\{y_m'\}$ equivalent if $\lim_\omega d(y_m,y_m')=0$. Denote by $[(y_m)]$ the equivalence class of $\{y_m\}$. The ultra-completion $Y_\omega$ of $Y$ with respect to $\omega$ is the metric space given by the set
$$Y_\omega:=\{[(y_m)]:\{y_m\} \text{ bounded sequence in }Y\},$$
equipped with the metric 
$$d_\omega([(y_m)],[(y_m')]):=\lim_\omega d(y_m,y_m').$$ 
The ultra-completion $Y_\omega$ of $Y$ is a complete metric space, even if $Y$ itself is not complete. 

Following~\cite{gw17}, a metric space $Y$ is said to be 1-complemented in some metric space $Z$ if $Y$ isometrically embeds into $Z$ and if there exists a 1-Lipschitz retraction from $Z$ to $Y$. By~\cite[Proposition 2.1]{gw17}, the class of metric spaces $Y$ which are 1-complemented in every ultra-completion of $Y$ includes proper metric spaces, NPC spaces, dual Banach spaces and injective metric spaces.

\subsection{Metric space valued Sobolev spaces}

\subsubsection{Korevaar-Schoen Sobolev spaces}

For each $\varepsilon>0$, we define an approximating energy $E_\varepsilon(u)\colon C_0(X)\to \R$ by
\begin{align*}
	E_\varepsilon(u)(f)=\int_{X}\dashint_{B(x,\varepsilon)}f(x)\frac{d(u(x),u(y))^2}{\varepsilon^2}d\mu(x)d\mu(y),
\end{align*}

\bd(Korevaar-Schoen-Sobolev space)\label{def:Korevaar-Schoen Sobolev}
Let $u\in L^2(X,Y)$. We say $u$ is in the Korevaar-Schoen-Sobolev space $KS^{1,2}(X,Y)$ if 
\begin{equation}\label{eq:Korevaar-Schoen energy}
	E(u):=\sup_{\varphi\in C_c(X,[0,1])}\limsup_{\varepsilon\to 0}E_\varepsilon(u)(\varphi)
\end{equation}
is finite. For $u\in KS^{1,2}(X,Y)$, $E(u)$ is the Korevaar-Schoen energy functional of $u$.
\ed

\subsubsection{Kuwae-Shioya Sobolev spaces}

As in the definition of Korevaar-Schoen Sobolev spaces, we first introduce the approximating energy.
\bd[Approximating energy]\label{def:approximating energy}
Fix a positive function $f\in C_0(X)$ and an admissible rate function $b$. For each $u\in L^p(X,Y)$, we set
\begin{align*}
	E_{\varepsilon}^{b}(u)(f):=\frac{1}{2b(\varepsilon)}\int_Xf(x)\int_{B(x,\varepsilon)}\frac{d(u(x),u(y))^2}{\varepsilon^2}d\mu(y)d\mu(x). 
\end{align*}
\ed 

If $X$ satisfies the SMCPBG, then it follows from~\cite[Theorem 3.1]{ks03} that for each $f\in C_0(X)$, the limit
\begin{align*}
	E^b(u)(f):=\lim_{\varepsilon\to 0}E^b_\varepsilon(u)(f)
\end{align*} 
exists. The limit functional $E^b(u)$ is called the Kuwae-Shioya energy functional of $u$.

\bd[Kuwae-Shioya energy]\label{def:energy functional of Kuwae-Shioya}
The Kuwae-Shioya energy of $u$ is defined to be 
$$E^b(u):=\sup_{0\leq f\leq 1, f\in C_0(X)}E^b(u)(f).$$
\ed 

The Sobolev space of Kuwae-Shioya is then defined to be 
\begin{align*}
	W^{1,2}(X,Y)=\Big\{u\in L^2(X,Y): E^b(u)<\infty \Big\}.
\end{align*}

\subsubsection{Sobolev spaces based on upper gradients}\label{subsec:Sobolev spaces}
Let $X=(X,d,\mu)$ be a metric measure space. Let $\Gamma$ a family of curves in $X$.  A Borel function $\rho\colon X\rightarrow [0,\infty]$ is admissible for $\Gamma$ if for every locally rectifiable curves $\gamma\in \Gamma$,
\begin{equation}\label{admissibility}
	\int_\gamma \rho\,ds\geq 1\text{.}
\end{equation}
The $2$-modulus of $\Gamma$ is defined as
\begin{equation*}
	\modulus_2(\Gamma) = \inf_{\rho} \left\{ \int_X \rho^2\,d\mu:\text{$\rho$ is admissible for $\Gamma$} \right\}.
\end{equation*}
A family of curves is called $2$-exceptional if it has $2$-modulus zero. We say that a property of curves holds for $2$-almost every curve if the collection of curves for which the property fails to hold is $2$-exceptional.

Let $X=(X,d,\mu)$ be a metric measure space and $Z=(Z,d_Z)$ be a metric space. A Borel function $g\colon X\rightarrow [0,\infty]$ is called an upper gradient for a map $u\colon X\to Z$ if for every rectifiable curve $\gamma\colon [a,b]\to X$, we have the inequality
\begin{equation}\label{ugdefeq}
	\int_\gamma g\,ds\geq d_Z(u(\gamma(b)),u(\gamma(a)))\text{.}
\end{equation}
If inequality \eqref{ugdefeq} holds for $2$-almost every curve, then $g$ is called a $2$-weak upper gradient for $u$.  

The concept of upper gradient was first introduced in~\cite{hk98} and then functions with $2$-integrable upper gradients were studied in~\cite{km98}. Later, the theory of real-valued Sobolev spaces based on upper gradients was explored in-depth in~\cite{s00}.

A $2$-weak upper gradient $g$ of $u$ is minimal if for every $2$-weak upper gradient $\tilde{g}$ of $u$, $\tilde{g}\geq g$ $\mu$-almost everywhere.  If $u$ has an upper gradient in $L^2_{\loc}(X)$, then $u$ has a unique (up to sets of $\mu$-measure zero) minimal $2$-weak upper gradient.  We denote the minimal upper gradient by $g_u$. The Newtonnian Sobolev space $N^{1,2}(X)$ consists of all functions $u\in L^2(X)$ with a minimal 2-weak upper gradient $g_u\in L^2(X)$. The notation $N^{1,2}_0(X)$ represents the class of all functions $g\in N^{1,2}(X)$ with compact support.

\bd[Sobolev capacity]\label{def:Sobolev capacity}
The 2-capacity of a set $E\subset X$ is defined to be the (possibly infinite) number
\begin{equation}\label{eq:def for 2 capacity}
	\capacity_2(E)=\inf \Big(\int_X |u|^2+g_u^2d\mu \Big),
\end{equation}
where the infimum is taken over all functions $u\in N^{1,p}(X)$ such that $u\geq 1$ outside a $2$-exceptional set of measure zero.
\ed 

Each function $u\in N^{1,2}(X)$ is quasicontinuous, meaning that for every $\varepsilon>0$, there exists an open set $G_\varepsilon\subset X$ with $\capacity_2(G_\varepsilon)<\varepsilon$ such that the restriction $u|_{X\backslash G_\varepsilon}$ is continuous.
Let $u,v\colon X\to Y$ be two mappings. We say that $u=v$ quasi-everywhere or q.e. in $X$ if $\capacity_2(\{x\in X:u(x)\neq v(x)\})=0$.

We say that a metric measure space $X=(X,d,\mu)$ supports a weak (1,2)-Poincar\'e inequality if there exist constants $C\geq 1$ and $\tau\geq 1$ such that
\begin{equation}\label{eq:Poin ineq}
	\dashint_B|u-u_B|d\mu\leq C\diam(B)\Big(\dashint_{\tau B}g^2d\mu\Big)^{1/2}
\end{equation}
for all open balls $B$ in $X$, for every function $u\colon X\to\bR$ that is integrable on balls and for every upper gradient $g$ of $u$ in $X$.

The Sobolev capacity enjoys many nice properties, such as the Cartan, Choquet and Kellogg properties, making it very useful in the theory of fine topology in Euclidean spaces; see for instance \cite[Section 2.4]{mz97}. These nice properties of Sobolev 2-capactiy remain valid in complete metric spaces that are doubling and support a weak (1,2)-Poincar\'e inequality; see e.g. \cite[Section 7]{bbl18} and \cite[Chapter 11]{bb11}. 

As an immdeiate consequence of the Choquet property of 2-capacity, one obtains that finely open sets are quasiopen. Here, a set $U\subset X$ is quasiopen if for every $\varepsilon>0$, there is an open set $G\subset X$ such that $\capacity_2(G)<\varepsilon$ and $G\cup U$ is open; see \cite{bbl18} for the definition of finely open sets and fine topology, and proofs of the above fact in the metric space setting. The aforementioned fact implies the following lemma, which substitutes for the partition of unity in metric spaces. This lemma was first proved (in a slightly stronger form) in \cite[Lemma 2.4]{km92} in Euclidean spaces and the proof remain valid in metric spaces with the help of the above fact; see also \cite[Lemma 2.153]{mz97}.

\bl\label{lemma:Kilpelainen-Maly}
Let $X$ be a complete doubling metric space that supports a weak (1,2)-Poincar\'e inequality. Let $\{U_{\alpha}\}_{\alpha}$ be a covering of an open set $U\subset X$ with each $U_\alpha$ being quasiopen. For each positive $g\in N^{1,2}_0(U)$, there exists a sequence of positive functions $g_j$ converging to $g$ in $N^{1,2}_0(U)$ such that each $g_j$ is a finite sum of functions in $N^{1,2}_0(U_\alpha)$.
\el

By \cite[Theorem 11.40]{bb11}, every quasicontinuous function $g\colon U\to [-\infty,\infty]$ is finely continuous at q.e. $x\in U$, which means that $g$ is continuous when $U$ equipped with the fine topology and $\overline{R}$ with the usual toplogy. In particular, the level sets of each function $g\in N^{1,2}(X)$ are finely open and thus quasiopen if $X$ is complete doubling with a weak (1,2)-Poincar\'e inequality.

The definition of metric space-valued Newtonnian Sobolev spaces $N^{1,2}(X,Y)$ can be found in~\cite[Chapter 7.1]{hkst12} and we do not recall it here. For each $u\in N^{1,2}(X,Y)$, we shall use $E^g(u)$ to denote the upper gradient energy functional of $u$, that is, 
$$E^g(u)=\int_X g_u^2d\mu.$$

\subsubsection{Equivalence of different Sobolev spaces}
There are also other definitions of metric space valued Sobolev spaces. Under mild assumptions for the source space $X$, one can show that all these definitions of Sobolev spaces coincide; see \cite[Chapter 10]{hkst12}. For our purpose, we shall point out a few facts on the relation of Korevaar-Schoen and Kuwae-Shioya Sobolev spaces with Newtonnian Sobolev spaces.

The following result shows the connection between $KS^{1,2}(X,Y)$ and $N^{1,2}(X,Y)$; see~\cite[Theorem 10.4.3 and Corollary 10.4.6]{hkst12}.

\bt\label{thm:KS is included in Newton}
Assume that $X$ is doubling. Then each $u\in KS^{1,2}(X,Y)$ has a $\mu$-representative $\tilde{u}$ in $N^{1,2}(X,Y)$ satisfying 
\begin{align*}
	E^g(\tilde{u})\leq CE(u),
\end{align*}
where the constant $C$ depends only on the doubling constant of $\mu$. If in addition $X$ supports a weak $(1,2)$-Poincar\'e inequality, then each $u\in N^{1,2}(X,Y)$ belongs to $KS^{1,2}(X,Y)$. 
\et

The following lemma is a simple consequence of \cite[Theorem 4.2]{ks03} and \cite[Theorem 10.3.4]{hkst12}.
\bl\label{lemma:ks in Newtion}
If $X$ is a compact metric space satisfying the SMCPBG, then there exists a constant $C$, depending only on the doubling constant of $X$ and the constant associated to the weak Poincar\'e inequality, such that for each $u\in W^{1,2}(X,Y)$, we have
\begin{align*}
	E^b(u)\leq CE^g(u).
\end{align*}
In particular, $W^{1,2}(X,Y)\subset N^{1,2}(X,Y)$ (up to $\mu$-representative).
\el
\begin{proof}
	Since $X$ is compact, it follows from~\cite[Theorem 4.2]{ks03} that $X$ is doubling with constant $c_X$ and supports a weak $(2,2)$-Poincar\'e inequality with constant $C_X$. In particular, each $u\in W^{1,2}(X)$ belongs to $P^{1,2}(X)$, the Sobolev space defined via weak (1,2)-Poincar\'e inequality (see e.g.~\cite[Section 10.3]{hkst12}). The claim follows then from~\cite[Theorem 10.3.4]{hkst12}.
\end{proof}

%\bd\label{def:quasicontinuous function}
%A mapping $u\colon X\to Y$ is said to be 2-quasicontinuous or just quasicontinuous if, for each $\varepsilon>0$, there exists an open set $G_\varepsilon\subset X$ with $\capacity_2(G_\varepsilon)<\varepsilon$ such that the restriction $u|_{X\backslash G_\varepsilon}$ is continuous. 
%\ed 
%
%By~\cite[Theorem 7.4.2 and Theorem 8.2.1]{hkst12}, when $X$ is a complete PI space (we refer doubling metric spaces supporting a weak $(1,2)$-Poincar\'e inequality as PI-spaces), each $u\in N^{1,2}(X,Y)$ is quasicontinuous. Combining this fact with Theorem~\ref{thm:KS is included in Newton}, we infer that when $X$ is a complete PI-space, each $u\in KS^{1,2}(X,Y)$ can be adjusted on a set of $\mu$-measure zero to become quasicontinuous.

\subsection{Cheeger energy and $\RCD(K,N)$ spaces}
In the celebrated works \cite{lv09,s06}, the notion of a lower bound on the Ricci curvature was introduced for general metric measure spaces. These spaces were termed CD$(K,N)$ or CD$(K,\infty)$ spaces, where the word CD refers to curvature-dimension, as it was defined via optimal transportation. In another significant work \cite{ags14}, a rich theory of (Sobolev) calculus was developed on such spaces. In particular, for each $u\in N^{1,2}(X)$, the Cheeger energy is defined as
\begin{equation}\label{def:cheeger energy}
\Ch(u):=\inf\Big\{\liminf_{n\to \infty}\frac{1}{2}\int_X|Du_n|^2d\mu:u_n\in \lip_b(X), u_n\to u\text{ in }L^2(X) \Big\}.
\end{equation}
One weakness for such spaces is that the Cheeger energy functional is not necessarily a quadratic form. In \cite{ags14b}, the authors introduced the stronger notion of a lower bound on the Riemannian Ricci curvature on a general metric measure space. Such metric spaces were called $\RCD(K,N)$ or $\RCD(K,\infty)$ spaces. Roughly speaking, $\RCD(K,N)$ spaces are CD$(K,N)$ plus the requirement that the Cheeger energy function $\Ch$ is quadratic; see \cite{ags14b,ags15} for more information about these spaces.

%If $(X,d,\mu)$ is $\RCD(K,N)$, then there is a unique bilinear symmetric form $\mathcal{E}\colon \mathbb{D}(\mathcal{\varepsilon})\times \mathbb{D}(\mathcal{\varepsilon})\to \R$

\section{Existence and uniqueness of solutions to the Dirichlet problem of Korevaar and Schoen}\label{sec:existence and uniqueness for KS}

\subsection{Formulation of Dirichlet problem associated to an energy functional $\mathcal{E}$}
Fix a mapping $\phi$ in a Sobolev space $S^{1,2}(X,Y)$ with energy functional $\mathcal{E}$ and we set
\begin{align}\label{eq:admissible mapping class}
S^{1,2}_{\phi}(\Omega,Y)=\{u\in S^{1,2}(X,Y): \text{ a } \mu\text{-reprentative of } u \text{ equals }\phi \text{ q.e. in } X\backslash \Omega \}.
\end{align}
In practice, we shall choose $(S^{1,2}(X,Y),\mathcal{E})$ to be the Korevaar-Schoen Sobolev space $(KS^{1,2}(X,Y),E)$, the Kuwae-Shioya Sobolev space $(W^{1,2}(X,Y),E^b)$, or the Newtonnian Sobolev space $(N^{1,2}(X,Y),E^g)$. As was pointed out in Section \ref{sec:preliminaries}, these Sobolev spaces coincide and the corresponding energy functionals are comparable on compact subsets.
%\begin{equation}\label{eq:admissible mapping class}
%KS^{1,2}_{\phi}(\Omega,Y)=\{u\in KS^{1,2}(X,Y) \text{ quasicontinuous}: u=\phi \text{ quasi-everywhere in } X\backslash \Omega \}.
%\end{equation}

Recall that the space $N^{1,2}_0(\Omega)$ consists of functions $u\in N^{1,2}(X)$ such that $u=0$ quasi-everywhere in $X\backslash \Omega$. We may equivalently characterize mappings in $S^{1,2}_\phi(\Omega,Y)$ as those mappings $u\in S^{1,2}(X,Y)$ that has a $\mu$-representative $\hat{u}$ for which $d(\hat{u},\phi)\in N^{1,2}_0(\Omega)$. Indeed, if $u\in S^{1,2}_\phi(\Omega,Y)$, then a $\mu$-representative $\hat{u}$ of $u$ satisfies $d(\hat{u},\phi)\in N^{1,2}(X)$ and $d(\hat{u},\phi)=0$ quasi-everywhere in $X\backslash \Omega$, which implies that $d(\hat{u},\phi)\in N_0^{1,2}(\Omega)$. On the other hand, if $u\in S^{1,2}(X,Y)$ has a $\mu$-representative $\hat{u}$ such that $d(\hat{u},\phi)\in N_0^{1,2}(\Omega)$, then $\hat{u}=\phi$ quasi-everywhere in $X\backslash \Omega$ and so $u\in S^{1,2}_\phi(\Omega,Y)$.

We say that a domain $\Omega\subset X$ supports a (1,2)-Poincar\'e inequality if there exits a constant $C_\Omega>0$ such that 
\begin{equation}\label{eq:omega weak poincare}
	\int_{\Omega} |v(x)|d\mu(x)\leq C_\Omega \Big(\int_{\Omega}|g_v(x)|^2d\mu(x)\Big)^{1/2} 
\end{equation}
for all functions $v\in N_0^{1,2}(\Omega)$.

\subsection{Energy functionals with property $\B$}\label{subsec:metric spaces with property B}

\bd[Property $\B$]\label{def:property B}
We say that an energy functional $\mathcal{E}$ between a complete metric measure space $X=(X,d,\mu)$ and a metric space $Y$ has property $\B$ if the following conditions are satisfied:
\begin{itemize}
	\item[(B1):] The Sobolev space $S^{1,2}(X,Y)$ is equivalent with $N^{1,2}(X,Y)$ and there exists a constant $C\geq 1$ such that for each $u\in S^{1,2}(X,Y)$
	$$C^{-1}\mathcal{E}(u)\leq E^g(u)\leq C\mathcal{E}(u).$$
	
	\item[(B2):] The metric measure space $(X,d,\mu)$ is locally a PI-space: for each relatively compact domain $K\subset X$, there exist a constant $c_K\geq 1$ and a radius $R_K>0$ 
	\begin{align*}
		\mu(B(x,2r))\leq c_d\mu(B(x,r))
	\end{align*}
	for each open ball $B(x,r)\subset K$ with $r< R_K/2$. Moreover, $X$ supports a the following local $(1,2)$-Poincar\'e inequality, for each compact set $K\subset X$, there exist $C_K>0$ and $\lambda_K\geq 1$ such that,
	\begin{align}\label{eq:1-2 Poincare inequality}
	\dashint_{B}|u(x)-u_B|d\mu(x)\leq C_K\diam B\Big(\dashint_{\lambda_K B}g^2d\mu(x)\Big)^{1/2},
	\end{align}
	for all open balls $B$ in $K$ with $\lambda_kB\subset K$, for every function $u\colon X\to\bR$ that is integrable on balls and for every upper gradient $g$ of $u$ in $X$.

	\item[(B3):] The energy functional $\mathcal{E}$ is lower semicontinuous with respect to the $L^2$-convergence, that is, if $u_i\to u$ in $L^2(X,Y)$, where $u_i$ and $u$ belong to $S^{1,2}(X,Y)$, then 
	\begin{align*}
		\mathcal{E}(u)\leq \liminf_{i\to \infty}\mathcal{E}(u_i).
	\end{align*}
 	
\end{itemize}
\ed

We would like to point out that condition (B2) above, volume doubling and weak (1,2)-Poincar\'e inequality, were initially introduced in \cite{j97} as one of the basic conditions to build up a  H\"older regularity theory for harmonic mappings in the singular space setting.

For the Korevaar-Schoen energy functional $E$ and the Kuwae-Shioya energy functional $E^b$, we introduce a slightly stronger notion.
\bd[Strong property $\B$]\label{def:strong property B}
We say that the Korevaar-Schoen energy functional $E$ between a complete metric space $X$ and a metric space $Y$ has strong property $\B$ if $E$ has property $\B$ and if in addition for each $u\in KS^{1,2}(X,Y)$ and each $f\in C_0(X)$, the pointwise limit $E(u)(f)=\lim_{\varepsilon\to 0}E_\varepsilon(u)(f)$ exists. The definition for the Kuwae-Shioya energy functional $E^b$ is similar by replacing $E$ with $E^b$.	
\ed

For Lipschitz manifolds considered in~\cite{g98}, admissible Riemannian polyhedrons considered in~\cite{c95,ef01,dm10} and metric spaces with the SMCP considered in~\cite{s98}, the Korevaar-Schoen energy functional $E$ has strong property $\B$.

\begin{example}\label{ex:metric space B}
On the following metric measure spaces $X$, the Korevaar-Schoen energy functional $E$ between $X$ and any metric space $Y$ have strong property $\B$.
\begin{itemize}
	\item[1).] $X$ is an admissible Riemannian polyhedron.
	\item[2).] $X$ satisfies the SMCP.
	\item[3).] $X$ is an $n$-dimensional Lipschitz submanifold of $\R^N$.
	\item[4).] $X$ is $\RCD(K,N)$ in the sense of \cite{gt20}\footnote{By the same reasoning, one can show that the strongly rectifiable metric spaces introduced in \cite{gt20} has property $\B$. Since we do not know any example other than $\RCD(K,N)$ spaces that is strongly rectifiable, we shall not consider these spaces in this article.}.
%	\item[4).] $X$ is an $n$-dimensional Alexandrov space with curvature bounded from below by some $\kappa\in \R$.
\end{itemize}	
\end{example}

\begin{proof}
	1). On admissible Riemannian polyhedrons, $E$ has strong property $\B$ follows from~\cite[Section 9]{ef01}.
	
	2). By~\cite[Proposition 4.5 and Theorem 6.4]{s98}, $(X,d,\mu)$ is locally doubling and supports a local $(2,2)$-Poincar\'e inequality, and hence it satisfies (B1) and (B2). The lower semicontinuity, property (B3), is a consequence of~\cite[Proof of Theorem 3.2]{ks03}. The additionally required property follows from the proofs of Theorem 3.1 and Theorem 4.1 in~\cite{ks03} by using the sub-partition lemma~\cite[Lemma 5.2]{s98}.
	%\textbf{We will include the details in the appendix?} 
	
	3). Property (B1) and (B2) are immediate as being a PI-space and Property (B3) follows from~\cite[Theorem 2]{g98}. The additionally required property is a direct consequence of~\cite[Theorem 1]{g98}.  
	
%	4). The locally doubling property of the volume measure is a simple consequence of the volume comparison. By~\cite[Theorem 4.2]{ks03}, $(X,d,\mu)$ supports a local $(p,p)$-Poincar\'e inequality, and hence it satisfies (B1). Property (B2) follows directly from~\cite[Theorem 3.1 and Theorem 4.1]{ks03}. The lower semicontinuity, property (B3), is a consequence of~\cite[Theorem 3.2]{ks03}. 
    4). Property (B1) and (B2) follow from \cite{r12,s06} and Property (B3) follows from \cite[Theorem 4.16]{gt20}. The additionally required property follows from \cite[Theorem 3.13]{gt20}.
\end{proof}

Regarding the other energy functionals, we have the following results.

\begin{example}\label{ex:SMCPBG imply property B}
1). If $(X,d,\mu)$ is a compact metric space satisfying the SMCPBG, then the Kuwae-Shioya energy functional $E^b$ between $X$ and any metric space $Y$ has strong property $\mathcal{B}$.

2). If $(X,d,\mu)$ is a compact PI space, then the upper gradient energy functional $E^g$ betwen $X$ and any metric space $Y$ has property $\mathcal{B}$.
\end{example}

\begin{proof}
1). Note that property (B1) and (B2) follow from~\cite[Theorem 4.2]{ks03} and the lower semicontinuity, property (B3), is a consequence of~\cite[Theorem 3.2]{ks03}. The additionally required property is a direct consequence of \cite[Theorem 3.1 and Theorem 4.1]{ks03}.

2). Properties (B1) and (B2) are clear. The lower semicontinuity follows from \cite[Theorem 7.3.9]{hkst12}. 
\end{proof}

%\begin{example}\label{ex:UG PI imply property B}
%If $(X,d,\mu)$ is a compact PI space, then the upper gradient energy functional $E^g$ satisfies property $\mathcal{B}$.
%\end{example}

\br\label{rmk:on metric subpartition to property B}
We would like to point out that it is very difficult, in general, to determine whether  an energy functional $\mathcal{E}$ between a given metric measure space $X$ and a metric space $Y$ has property $\B$ or not and often it requires a very nice geometric structure on $X$ (but not necessarily for $Y$) in order to admit such a nice energy functional.
%More generally, if a metric space $X$ is sufficiently nice so that a general version of the sub-partition lemma as in~\cite[Lemma 3.2]{ks03} holds  (see also~\cite[Lemma 1.3.1]{ks93} and~\cite[Lemmas 5.2 and 5.3]{s98} for simpler versions), then $X$ satisfies properties (B2) and (B3). Indeed, one can directly verify that the proofs of Theorem 3.1 and Theorem 3.2 in \cite{ks03} work as long as the sub-partition lemma holds.
\er

\subsection{Proof of Theorem~\ref{thm:main thm existence and uniqueness}}

We will need the following compactness result from~\cite[Theorem 3.1]{gw17} in our existence proof below.
\bt[Generalized Rellich compactness]\label{thm:Rellich compactness}
Suppose $\mathcal{E}$ is an energy functional between a compact metric measure space $X$ and a metric space $Y$ that has property $\B$. Let $\{u_m\}\subset KS^{1,2}(X,Y)$ be a sequence such that 
\begin{align}\label{eq:generalized compactness}
	\sup_{k\in \mathbb{N}}\int_{X}d^2(u_k(x),y_0)d\mu(x)+\mathcal{E}(u_k)<\infty 
\end{align}
for some $y_0\in Y$. Then after possibly passing to a subsequence, there exist a complete metric space $Z$, isometric embeddings $\varphi_k\colon Y\to Z$, and $v\in S^{1,2}(X,Z)$ such that $\varphi_k\circ u_k$ converges to $v$ in $L^2(X,Z)$.
\et

\br\label{rmk:on generalized compactness} 
Note that in~\cite[Theorem 3.1]{gw17}, $X$ is assumed to be a bounded Lipschitz domain in $\R^n$ and $\mathcal{E}=E$ is the Korevaar-Schoen energy functional. However, this assumption was only used to ensure that each $u\in KS^{1,2}(X,Y)$ has a $\mu$-representative $\hat{u}$ that satisfies the pointwise inequality
\begin{equation}\label{eq:hajlasz estimate}
d(\hat{u}(x),\hat{u}(x'))\leq d(x,x')(h(x)+h(x'))
\end{equation}
for some $h\in L^2(X)$. In our setting, this fact is well-known by property (B1) and (B2) (see e.g.~\cite[Theorem 10.5.2 and Theorem 10.5.3]{hkst12}). 
\er

\begin{proof}[Proof of Theorem~\ref{thm:main thm existence and uniqueness}]
The proof follows closely the approach of~\cite[Proof of Theorem 1.4]{gw17}. Let $\{u_k\}\subset S^{1,2}_\phi(\Omega,Y)$ be an energy minimizing sequence. 
Up to a $\mu$-representative, we may assume that each $h_k(x)=d(u_k(x),\phi(x))\in N^{1,2}_0(\Omega)$. Since $\sup_{k}\mathcal{E}(h_k)<\infty$, it follows easily from property (B1) and the (1,2)-Poincar\'e inequality~\eqref{eq:Poin ineq} that $\sup_{k}\|h_k\|_{L^2(X)}<\infty$. Hence
\begin{align*}
\sup_{k}\int_{X} d^2(u_k(x),y_0)d\mu(x)+\mathcal{E}(u_k)<\infty.
\end{align*}

We may apply Theorem~\ref{thm:Rellich compactness} to find, after possibly passing to a subsequence, a complete metric space $Z=(Z,d_Z)$, 
isometric embeddings $\varphi_k\colon Y\to Z$ and $v\in S^{1,2}(X,Z)=KS^{1,2}(X,Z)$ such that 
%$\varphi_k(C_m)\subset Z^m$ for all $k,m$ and 
$v_k:=\varphi_k\circ u_k$ converges in $L^2(X,Z)$ to $v$ as $k\to \infty$. 
After passing to a further subsequence, we may assume that $v_m$ converges almost everywhere on $X$. 
Let $N\subset X$ be a set of $\mu$-measure zero such that $v_k(z)\to v(z)$ for all $z\in X\backslash N$.

By our assumption on $Y$, there exists an ultra-completion $Y_\omega$ on $Y$ such that $Y$ admits a 1-Lipschitz retraction $P\colon Y_\omega\to Y$. Define a subset of $Z$ by $B:=\{v(z):z\in X\backslash N \}$.
%$\cup \varphi(C)$. 
The map $\psi\colon B\to Y_\omega$, given by $\psi(v(z))=[(u_k(z))]$ when $z\in X\backslash N$
% and $\psi(\varphi(y))=[(y)]$ when $y\in C$ 
is well-defined and isometric by~\cite[Lemma 2.2]{gw17}. Since $Y_\omega$ is complete, there exists a unique extension of $\psi$ to $\overline{B}$, which we denote again by $\psi$. After possibly redefining the map $v$ on $N$, we may assume that $v$ has image in $\overline{B}$ and hence $v$ is an element of $KS^{1,2}(X,\overline{B})$. Now, we define a mapping by $u:=P\circ \psi\circ v$ and then $u$ belongs to $KS^{1,2}(X,Y)=S^{1,2}(X,Y)$ and satisfies
\begin{align*}
E(u)\leq E(v)\leq \lim_{k\to \infty}E(v_m)=\lim_{k\to \infty} E(u_m)
\end{align*}
by the lower semicontinuity of the energy from property (B2). 

It remains to show that $u=\phi$ quasi-everywhere in $X\backslash \Omega$. For each $m\in \mathbb{N}$, we write $B_m=\{v_m(z)\colon z\in X\backslash N\}$ and define $\psi_m\colon B_m\to Y_\omega$ by $\psi_m(v_m(z))=[(u_m(z))]$ (the constant sequence). Then $\psi_m$ is well-defined and isometric by~\cite[Lemma 2.2]{gw17}. Moreover, $u_m=P\circ \psi_m\circ v_m$. Note that for q.e. $x\in X\backslash \Omega$, $u_m(x)=\phi(x)$ for all $m\in \mathbb{N}$. 
So for all such point $x$ we have
\begin{align*}
	d(\phi(x),u(x))&=d(u_m(x),u(x))=d(P\circ \psi_m\circ v_m(x),P\circ \psi\circ v(x))\\
	&\leq d(\psi_m\circ v_m(x),\psi\circ v(x))=\lim_{k\to \omega}d(u_m(x),u_k(x))\\
	&=d(\phi(x),\phi(x))=0,
\end{align*}
which implies that $u(x)=\phi(x)$ q.e. in $X\backslash \Omega$. Thus $u\in S^{1,2}_{\phi}(X,Y)=KS^{1,2}_\phi(X,Y)$ as required.
\end{proof}

\subsection{Proof of Theorem \ref{thm:uniqueness}}

The following lemma will be the key for the proof of Theorem~\ref{thm:uniqueness}. 

\bl\label{lemma:key lemma for existence}
Let $X$ be a complete locally PI-space and $Y$ a metric space. Let $\{u_k\}$ be a sequence in $KS^{1,2}_\phi(X,Y)$ with uniformly bounded energy. Suppose $u_k$ converges in $L^2(X,Y)$ to some $u\in KS^{1,2}(X,Y)$. Then $u\in KS^{1,2}_\phi(X,Y)$.
\el 
\begin{proof}
	We only need to show that a $\mu$-representative of $u$ equals $\phi$ quasi-everywhere in $X\backslash \Omega$, or equivalently, for each compact set $K\subset X$, $u=\phi$ q.e. in $K$. Thus we may assume without loss of generality that $X$ is a complete PI space. By embedding $Y$ isometrically into $l^\infty(Y)$ if necessary, we may further assume that $Y$ is (isometrically) contained in a Banach space. 
	
	Note that Theorem~\ref{thm:KS is included in Newton} implies that $g_{u_k}$ is bounded in $L^2(X)$. Taking a further subsequence if necessary, we may assume that $g_{u_k}\rightharpoonup g$ weakly in $L^2(X)$. By Mazur's lemma (see e.g.~\cite[Section 2.3]{hkst12}), a convex combination $\hat{g}_i=\sum_{k=i}^{N_i} a_{ki}g_{u_k}$ of $g_{u_k}$ converges to $g$ strongly in $L^2(X)$. Consequently, the sequence $\hat{u}_i=\sum_{k=i}^{N_i} a_{ki}u_k$ converges to $u$ in $L^2(X,Y)$ with $\hat{g}_i$ being a 2-weak upper gradient of $\hat{u}_i$. Then~\cite[Proposition 7.3.7]{hkst12} implies that $u$ has a $\mu$-representative in $N^{1,2}(X,Y)$ with each Borel representative of $g$ as a 2-weak upper gradient. Moreover, a subsequence of $\hat{u}_i$ converges pointwise to this representative of $u$ outside a set of 2-capacity zero. Note that in $X\backslash \Omega$, $u_i=\phi$ quasi-everywhere and so is $\hat{u}_i$. This implies that this $\mu$-representative of $u$ coincides with $\phi$ quasi-everywhere in $X\backslash \Omega$. In particular, $u\in KS^{1,2}_\phi(\Omega,Y)$. 
\end{proof}

\begin{proof}[Proof of Theorem~\ref{thm:uniqueness}]
We will follow the approach of~\cite{ks93} to show the existence and uniqueness of energy minimizers.

\textbf{Uniqueness:}
For two mappings $u,v\in KS^{1,2}(X,Y)$, we denote by $w$ the middle point mapping of $u$ and $v$. More precisely, for each $x\in X$, we set $w(x)$ to be the middle point of the geodesic connection $u(x)$ and $v(x)$. Note that $d(u,v)\in KS^{1,2}(X)$. We next show that $w\in KS^{1,2}(X,Y)$.

For $x,y\in X$, by~\cite[Equation (2.2iii)]{ks93},
\begin{align*}
	2d^2(w(x),w(y))&\leq d^2(u(x),u(y))+d^2(v(x),v(y))\\
	&\qquad\qquad\qquad-\frac{1}{2}[d(u(y),v(y))-d(u(x),v(x))]^2.\numberthis\label{eq:2.2iii}
\end{align*}
Integrating and averaging~\eqref{eq:2.2iii} on the ball $B(x,\varepsilon)$ with respect to $y$; and then multiplying by $f(x)$, where $f\geq 0$ and $f\in C_0(\Omega)$, then integrating with respect to $x$ and sending $\varepsilon$ to zero, we deduce from Property (B3) that 
\begin{align*}
	2E(w)(f)\leq E(u)(f)+E(v)(f)-\frac{1}{2}E(d(u,v))(f),
\end{align*}
or equivalently,
\begin{align}\label{eq:3.2}
	2E(w)\leq E(u)+E(v)-\frac{1}{2}E(d(u,v))
\end{align}
from which we infer $w\in KS^{1,2}(X,Y)$.

Suppose $u$ and $v$ are two solutions of the Dirichlet problem. Then the middle point mapping $w$ of $u$ and $v$ belongs to $KS^{1,2}_\phi(\Omega,Y)$ since $w=\phi$ quasi-everywhere in $X\backslash \Omega$. It follows that $E(w)\geq E(u)=E(v)$. By~\eqref{eq:3.2}, we infer that $E(h)=0$, where $h=d(u,v)$, and so for each relatively compact domain $K$ of $X$ that contains $\Omega$, we have by property (B1) and Theorem~\ref{thm:KS is included in Newton} that
\begin{align*}
	\int_{K}g_{h}^2d\mu\leq CE(h)=0. 
\end{align*}
Since $h=0$ quasi-everywhere on $X\backslash \Omega$, $h=0$ quasi-everywhere on $X\backslash K$ as well and so $g_h=0$ $\mu$-a.e. on $X\backslash K$. We thus conclude that
 \begin{align*}
 \int_{X}g_{h}^2d\mu\leq CE(h)=0. 
 \end{align*}
In particular, $g_h=0$ $\mu$-a.e. in $X$. Since $X$ is locally a PI-space, we conclude that $h=c$ for some constant $c\in \R$ $\mu$-a.e. in $X$. Since $h=0$ quasi-everywhere in $X\backslash \Omega$, we conclude that $c=0$. Consequently, $u=v$ $\mu$-a.e. in $X$. This shows the uniqueness of solutions for the Dirichlet problem.

\textbf{Existence:} Let $\{u_i\}$ be an energy minimizing sequence for the Dirichlet problem, i.e.,
\begin{align*}
\lim_{i\to \infty} E(u_i)=\inf_{v\in KS^{1,2}_\phi(\Omega,Y)}E(v)=:E_0.
\end{align*}
%We may assume that each $u_i$ is quasicontinuous. 
Then $u_{ij}:=d(u_i,u_j)\in N^{1,2}_0(\Omega)$ and the middle point mapping $w_{ij}$ of $u_i$ and $u_j$ belongs to $KS^{1,2}_\phi(\Omega,Y)$. Thus
\begin{align*}
E(u_i)+E(u_j)-2E(w_{ij})\leq E(u_i)+E(u_j)-2E_0.
\end{align*}
By~\eqref{eq:3.2}, this implies that
\begin{align}\label{eq:3.3}
\int_{X}g_{u_{ij}}^2d\mu=\int_{\Omega}g_{u_{ij}}^2d\mu \leq CE(u_{ij})\to 0 \text{ as }i,j\to \infty.
\end{align}
Note that $g_{u_{ij}}=0$ in $X\backslash \Omega$, since $u_{ij}=0$ quasi-everywhere in $X\backslash \Omega$. Since $\Omega$ supports the weak $(1,2)$-Poincar\'e inequality~\eqref{eq:omega weak poincare}, we infer that $u_{ij}\to 0$ in $L^2(\Omega)$ and so also in $L^2(X,Y)$ as $i,j\to \infty$, and hence $\{u_i\}$ has a limit $u$ in the complete metric space $L^2(X,Y)$. By the lower semicontinuity of the energy functional, 
\begin{align*}
E(u)\leq \liminf_{i\to \infty} E(u_i)=E_0
\end{align*}
and so $u\in KS^{1,2}(X,Y)$. By Lemma~\ref{lemma:key lemma for existence}, $u\in KS^{1,2}_\phi(X,Y)$ and so $u$ is the required energy minimizer.

\end{proof}

\br\label{rmk:on existence and uniqueness}
1). The proof of Theorem~\ref{thm:uniqueness} above implies that the requirement that $\Omega$ supports a $(1,2)$-Poincar\'e inequality can be dropped if in additional $X$ supports a $(1,2)$-Poincar\'e inequality, i.e., there exists a constant $C_X>0$ such that
\begin{equation}\label{eq:X weak poincare}
\int_{X}|v(x)|d\mu(x)\leq C_X \Big(\int_{X}|g_v(x)|^2d\mu(x)\Big)^{1/2} 
\end{equation}
for all $v\in N^{1,2}_0(\Omega)$.

2). We can modify the proof of Theorem~\ref{thm:uniqueness}, similar as that in~\cite[Proof of Theorem 2 (a)]{f05}, so that it works for $Y$ being a complete metric space with curvature bounded from above (the so-called CAT$(k)$-spaces). But then one has to consider mappings into a closed geodesic ball $B$ in $Y$ with radius $R<\frac{\pi}{2\sqrt{k}}$. 
%Similar observations have been made by Jost~\cite[Theorem 5.2]{j96}.  

3). The proof shows that Theorem \ref{thm:uniqueness} holds for the Kuwae-Shioya energy functional if instead we assume $X$ satisfies the SMCPBG. 
\er

\section{H\"older Regularity of harmonic mappings}\label{sec:holder continuity of KS}

In this section, we study the interior regularity of solutions for the Dirichlet problem of Korevaar-Schoen, for which, we term harmonic mappings. We shall prove that harmonic mappings are locally H\"older continuous, provided the metric space $X$ attains certain analytic property and $Y$ is NPC. The proof of Theorem~\ref{thm:main thm regularity} 1) is based on a combination of the arguments of Jost~\cite{j97} and Lin~\cite{lin97}, while the proof of Theorem~\ref{thm:main thm regularity} 2) follows closely the approach of Jost~\cite{j97}.

We will need the theory of Dirichlet forms on Hilbert spaces to separate a class of metric spaces, which we name them as metric spaces with property $\Cp$.

\subsection{Dirichlet forms}

Recall that a Dirichlet form $\mathcal{E}$ on $L^2(X,\mu)$ is a closed nonnegative definite and symmetric bilinear form defined on a dense linear subspace $\mathbb{D}=\mathbb{D}(\mathcal{E})$ of $L^2(X,\mu)$, that satisfies the Markovian property
$$\mathcal{E}(v,v)\leq \mathcal{E}(u,u)\quad \text{for all }u\in \mathbb{D},$$
where $v=\min\{1,\max\{u,0\} \}$. A Dirichlet form $\mathcal{E}$ on $L^2(X,\mu)$ is said to be strongly local if $\mathcal{E}(u,v)=0$ whenever $u,v\in \mathbb{D}$ with $u$ a constant on a neighborhood of the support of $v$; to be regular if there exists a subset of $D\cap C_0(X)$ which is both dense in $C_0(X)$ with the uniform norm and in $\mathbb{D}$ with the graph norm $\|\cdot\|_{\mathbb{D}}$ defined by $\|u\|_{\mathbb{D}}=\sqrt{\int_X u^2d\mu+\mathcal{E}(u,u)}$ for each $u\in\mathbb{D}$.

By the construction of Beurling and Deny~\cite{fot94}, each regular strongly local Dirichlet form $\mathcal{E}$ on $L^2(X,\mu)$ can be written as 
\begin{equation}\label{eq:representation of Dirichlet form}
\mathcal{E}(u,v)=\int_X d\Gamma(u,v) \quad \text{for all }u,v\in \mathbb{D},
\end{equation}
where $\Gamma$ is an $\mathcal{M}(X)$-valued nonnegative definite and symmetric bilinear form
defined by the formula
\begin{equation}\label{eq:energy measure of Dirichlet form}
\int_X \varphi d\Gamma(u,v)=\frac{1}{2}\big[\mathcal{E}(u,\varphi v)+\mathcal{E}(v,\varphi u)-\mathcal{E}(uv,\varphi) \big]
\end{equation}
for all $u,v\in \mathbb{D}\cap L^\infty(X,\mu)$ and $\varphi\in \mathbb{D}\cap C_0(X)$. We call $\Gamma(u,u)$ the Dirichlet energy measure (squared gradient) and $\sqrt{\frac{d\Gamma(u,u)}{d\mu}(x)}$ the length of the gradient of $u$ at $x$.

For a strongly local Dirichlet form $\mathcal{E}$, its energy measure $\Gamma$ is local and satisfies the Leibniz rule and the chain rule. Both $\mathcal{E}(u,v)$ and $\Gamma(u,v)$ can be defined for $u,v\in \mathbb{D}_{\loc}(X)$, the collection of all $u\in L^2_{\loc}(X)$ satisfying that for each relatively compact set $K\subset X$, there exists a function $w\in \mathbb{D}$ such that $u=w$ almost everywhere on $K$. With this, the intrinsic distance on $X$ associated to $\mathcal{E}$ is defined by
\begin{align}\label{eq:intrinsic distance}
d_{\mathcal{E}}(x,y)=\sup\{u(x)-u(y):u\in\mathbb{D}_{\loc}(X)\cap C(X), \Gamma(u,u)\leq \mu \}.
\end{align}
Here $\Gamma(u,u)\leq \mu$ means that $\Gamma(u,u)$ is absolutely continuous with respect to $\mu$ and $\frac{d\Gamma(u,u)}{d\mu}\leq 1$ almost everywhere.

Given a Dirichlet form $\mathcal{E}$ on the Hilbert space $L^2(X,\mu)$, there exists a unique self-adjoint operator $A$ on $L^2(X,\mu)$ with the properties $\mathbb{D}=D(A^{1/2})$ and 
$$-\mathcal{E}(u,v)=(u,Av)=\int_X u(x)Av(x)d\mu(x)$$
for all $u\in \mathbb{D}$ and $v\in D(A)$.

\subsection{Approximating energy densities with Property $\Cp$}\label{subsec:metric spaces C}
\bd[Approximating energy density]\label{def:approximate energy density}
Fix an admissible rate function $c\colon X\times (0,\infty)\to (0,\infty)$. For each $\varepsilon>0$, we define an approximating energy density $e_\varepsilon^{c,Y}$ acting on $u\colon X\to Y$ as follows:
$$e_\varepsilon^{c,Y}(u)(x)=\frac{1}{c(x,\varepsilon)}\int_{B(x,\varepsilon)}\frac{d^2(u(x),u(y))}{\varepsilon^2}d\mu(y).$$
We simply write $e_\varepsilon^c$ when $(Y,d)=(\R,|\cdot|)$.
\ed

For each $u\in \lip_0(X)$, the space of Lipschitz functions on $X$ with compact support, we introduce the approximating energy $E_\varepsilon$ as 
\begin{align*}
E_\varepsilon^c(u)=\int_{X}e_\varepsilon^c(u)(x)d\mu(x).
\end{align*}
%\begin{align*}
%	E_\varepsilon(u)=\int_{X}\dashint_{B(x,\varepsilon)}\frac{|u(x)-u(y)|^2}{\varepsilon^2}d\mu(y)d\mu(x).
%\end{align*}

\bd[Property $\mathcal{C}$]\label{def:property C}
We say that the approximating energy density $e_\varepsilon^c$ on a complete metric space $X$ has property $\Cp$ if the following conditions are satisfied:
\begin{itemize}
	\item[(C1):] The pointwise limit $E_0=\lim_{\varepsilon\to 0}E_\varepsilon^c$ exists and induces a regular strongly local Dirichlet forms $\mathcal{E}_0$ on $L^2(X,\mu)$ via the formula 
	\begin{equation}\label{eq:def to generate Dirichlet form}
	\mathcal{E}_0(u,v)=\lim_{\varepsilon\to 0}\int_{X}\frac{1}{c(x,\varepsilon)}\int_{B(x,\varepsilon)}\frac{[u(x)-u(y)][v(x)-v(y)]}{\varepsilon^2}d\mu(y)d\mu(x).
	\end{equation}
	Moreover, the intrinsic distance $d_0:=d_{\mathcal{E}_0}$ associated to $\mathcal{E}_0$ is locally bi-Lipschitz equivalent with the original distance $d$ on $X$. In particular, it induces the same topology as the underlying topology on $X$.	
	
	\item[(C2):] Equip with the intrinsic metric, the space $(X,d_0,\mu)$ becomes a complete locally PI space. That is $\mu$ is a locally doubling measure on $(X,d_0)$: for each compact set $K\subset (X,d_0)$, there exists a doubling constant $c_K\geq 1$ such that 
	$$\mu(2B)\leq c_K\mu(B)\quad \text{ for each ball }B \text{ with } 2B\subset K\subset (X,d_0),$$
	and that $(X,d_0,\mu)$ supports a local weak $(2,2)$-Poincar\'e inequality: for each compact set $K\subset (X,d_0)$, there exist $C_K>0$ and $\lambda_K\geq 1$ such that for each ball $B\subset (X,d_0)$ with $\lambda_K B\subset K$ and all $u\in \mathbb{D}(\mathcal{E}_0)$,
	\begin{align}\label{eq:2-2 Poincare inequality}
	\int_{B}|u(x)-u_B|^2d\mu(x)\leq C_K(\diam B)^2\int_{\lambda_K B}d\Gamma_0(u,u),
	\end{align}
	where $\Gamma_0$ is the Dirichlet energy measure corresponding to the Dirichlet form $\mathcal{E}_0$.	
\end{itemize}
\ed

\br\label{rmk:on domain of Dirichlet forms}
	Note that if the approximating energy density $e_\varepsilon^c$ on a complete metric space $X$ has property $\Cp$, then it is well-known that $\mathbb{D}(\mathcal{E}_0)=N^{1,2}(X)$ and $\mathbb{D}_{\loc}(\mathcal{E}_0)=N^{1,2}_{\loc}(X)$. Moreover, $\lip_0(X)$ is dense in $\mathbb{D}(\mathcal{E}_0)$; see e.g.~\cite[Theorem 2.2]{kz12}.
\er

\bd[Strong property $\Cp$]\label{def:strong property c}
We say that the approximating energy density $e_\varepsilon^{c,Y}$ between a complete metric space $X$ and an NPC space $Y$ has strong property $\Cp$ if  
\begin{itemize}
	\item $e_\varepsilon^{c}$ has property $\Cp$;
	\item for each $u\in KS^{1,2}(X,Y)$, there exists an energy measure $\mu_{u,Y}$ such that $e_\varepsilon^{c,Y}(u)d\mu\rightharpoonup \mu_{u,Y}$ weakly, that is, 
	%the following representation formula holds: for each $u\in KS^{1,2}(X,Y)$, there exists a energy density $e(u)\in L^1(X)$ such that 
	for each positive function $\eta\in C_0(X)$,
	\begin{align}\label{eq:representation formula}
		\limsup_{\varepsilon\to 0}\int \eta(x)e_\varepsilon^{c,Y}(u)(x)d\mu(x)=\int \eta(x)d\mu_{u,Y}(x);
	\end{align}
\end{itemize}
\ed
In particular, if  the approximating energy density $e_\varepsilon^{c,Y}$ has strong property $\Cp$, then we have $E(u)(f)=\int_X f d\mu_{u,Y}(x)$ for each $u\in KS^{1,2}(X,Y)$ and each $f\in C_0(X)$. As an application of Theorem \ref{thm:KS is included in Newton}, we infer that for each compact $K\subset X$, there exists $C_K>0$ such that for $u\in KS^{1,2}(X,Y)$, 
\begin{equation}\label{eq:strong property C 3}
	C_K^{-1}\int_Kg_ud\mu\leq \limsup_{\varepsilon\to 0}\int_K e_\varepsilon^{c,Y}(u)d\mu\leq C_K^{-1}\int_Kg_ud\mu.
\end{equation}

Fix an NPC space $Y$, for a large class of metric measure spaces $X$, the Korevaar-Schoen approximating energy density $e_\varepsilon^Y$, i.e. $e_\varepsilon^Y=e_\varepsilon^{c,Y}$ with $c(x,\varepsilon)=\mu(B(x,\varepsilon))$ in Definition \ref{def:approximate energy density}, between  $X$ and $Y$, has strong property $\Cp$. 
\begin{example}\label{ex:metric space C}
Fix an NPC space $Y$. Then on the following $X$, the Korevaar-Schoen approximating energy density $e_\varepsilon^Y$ between $X$ and $Y$ have strong property $\Cp$.
	\begin{itemize}
		\item[1).] $X$ is an admissible Riemannian polyhedron.
		\item[2).] $X$ satisfies the SMCP.
		\item[3).] $X$ is an $n$-dimensional Lipschitz submanifold of $\R^N$.
		\item[4).] $X$ is $\RCD(K,N)$ in the sense of \cite{gt20}.
	\end{itemize}	
\end{example}

\begin{proof}
	1). Suppose $M$ is an admissible Riemannian polyhedron. Property (C1) is a direct consequence of~\cite[Proposition 7.3 and Theorem 9.1]{ef01} and (C2) follows from~\cite[Corollary 4.1, Thoerem 5.1]{ef01}. That $e_\varepsilon^{Y}(u)\rightharpoonup e^{Y}(u)$ is contained in~\cite[Theorem 9.1]{ef01}.	
	
	2). Suppose $X$ satisfies SMCP. Property (C1) follows from~\cite[Theorem 3.3 and Theorem 6.8]{s98}. Property (C2) is a direct consequence of~\cite[Proposition 4.5 and Theorem 6.4]{s98} by noticing that the original distance $d$ on $X$ is locally comparable with the intrinsic distance $d_0$ induced by the Dirichlet form $\mathcal{E}_0$, which is consequence of~\cite[Proposition 6.6 and Proposition 6.7]{s98}. That $e_\varepsilon^Y(u)\rightharpoonup e^Y(u)$ follows from~\cite[The proof of Theorem 6.1]{s98} by noticing that only a sub-partition lemma is needed for the proof (see also~\cite[Theorem 4.1]{ks03}).

	3). The first assertion in property (C1) follows from the subpartition lemma~\cite[Lemma 3]{g98} and the proof of Corollary 5.5 in~\cite{s98} or alternatively, the proof of Theorem 3.1 in~\cite{ks03}. Property (C2) is clear as being a Lipschitz manifold.

	The second claim in property (C1) follows by a minor modification of the proofs of Proposition 6.6 and Proposition 6.7 in~\cite{s98}, where one uses the local doubling property of the measure, the local weak (2,2)-Poincar\'e inequality and~\cite[Theorem 1]{g98}. Alternatively, one can follow the exact proof of Example \ref{exam:Kuwae-Shioya strong property c} below.
	
	That $e_\varepsilon(u)\rightharpoonup e(u)$ follows from~\cite[Theorem 1]{ef01}.
	
	%4). In this case, property (C1) is a consequence of~\cite[Theorem 3.1]{ks03}, where the local equivalence of the distance $d_0$ and $d$ follows from the proof of Lemma~\ref{lemma:SMCPBG imply property C} below, or essentially, the proof of Proposition 6.6 and Proposition 6.7 in~\cite{s98}. Property (C2) is a direct consequence of \cite[Theorem 4.2]{ks03}. That $e_\varepsilon(u)\rightharpoonup d\mu_u$ is contained in~\cite[Theorem 4.1]{ks03}.
	
	4). By \cite[Theorem 3.13]{gt20}, the pointwise limt $E_0=\lim_{\varepsilon}E_\varepsilon$ exists and by \cite[Proposition 4.19 and Theorem 4.1]{gt20}, for each $u\in \lip_0(X)$, we have
	$$E_0(u)=\int_X e_2(u)^2d\mu=c(d)\int_X|du|^2d\mu=c(d)\int_X|g_u|^2d\mu=c(d)\Ch(u),$$
	where $c(d)$ is a constant depending only on the dimension of $X$ and $\Ch(u)$ denotes the Cheeger energy of $u$ introduced in \cite{ags14} (see also \cite{ags15}). It is well-known that the Cheeger energy function $\Ch$ generates a strongly local regular Dirichlet form $\mathcal{E}_0$ by the following formula  
	\begin{equation}\label{eq:cheeger energy}
		\mathcal{E}_0(u,u)=2\Ch(u,u);
	\end{equation}
	see \cite[Section 4.3]{ags14b} or \cite{ags15}. By \cite[Theorem 6.10]{ags14b} the intrinsic distance $d_0:=d_{\mathcal{E}_0}$  is locally bi-Lipschitz equivalent with the original distance $d$ on $X$. This verifies property (C1).
	
	Since $d_0$ is locally bi-Lipschitz equivalent to the original metric and $\mu$ is locally doubling on $(X,d)$, it is also locally doubling on $(X,d_0)$. To prove the weak $(2,2)$-Poincar\'e inequality, first note that $(X,d,\mu)$ supports a weak $(1,2)$-Poincar\'e inequality (see \cite{r12}), and thus by \cite[Proposition 2.1 and Proposition 2.2]{kz12}, $(X,d_0,\mu)$ supports a weak $(1,2)$-Poincar\'e inequality as well. Since the space $(X,d_0,\mu)$ is locally quasiconvex and doubling, by \cite[Remark 9.1.19]{hkst12}, the weak $(2,2)$-Poincar\'e inequality \eqref{eq:2-2 Poincare inequality} follows from the weak $(1,2)$-Poincar\'e inequality. This verifies property (C2).
	
	Note that in our case, we actually have, by \cite[Theorem 3.13 and Proposition 4.6]{gt20}, for each $u\in KS^{1,2}(X,Y)$,
	$$e_{\varepsilon}(u)\to e_2(u)^2=\mathcal{S}_2(g_u)^2\qquad \text{in }L^2(X),$$
	where $\mathcal{S}_2$ is a seminorm on $\R^d$ defined as in \cite[Definition 3.12]{gt20} and $c(d)>0$ is a constant depending only on the dimension $d$ of $X$. This verifies the additionally required property.
	
\end{proof}

We next verify an example for the Kuwae-Shioya approximating energy density $e_\varepsilon^{b,Y}$, i.e. $e_\varepsilon^b=e_\varepsilon^{c,Y}$ with $c(x,\varepsilon)=b(\varepsilon)$ independent of $x$. 
\begin{example}\label{exam:Kuwae-Shioya strong property c}
Let $X$ be a complete metric measure space satisfying the SMCPBG and $Y$ an NPC space. Then the Kuwae-Shioya approximating energy density $e_\varepsilon^{b,Y}$ between $X$ and $Y$ has strong property $\Cp$.
\end{example}	

\begin{proof}
	The first assertion in property (C1) follows from~\cite[Theorem 3.1]{ks03}, where we write $E_0^b$ instead of $E_0$ and $\mathcal{E}_0^b$ instead of $\mathcal{E}_0$. 
	
	We next show the second assertion in property (C1), that is, the intrinsic distance $d_0$ associated to $\mathcal{E}_0^b$ satisfies $d_0\approx d$ on each compact set $K\subset X$. In particular, it induces the same topology as the underlying topology on $X$. The proof is similar to~\cite[Proofs of Propositions 6.5 and 6.6]{s98}.

	\underline{Claim 1:} $d_0\geq d$ on $X\times X$.
	
	Indeed, for each $x,z\in X$, each $y\in B(x,\varepsilon)$, the triangle's inequality implies that
	\begin{align*}
		\Big|\frac{d(x,z)-d(y,z)}{\varepsilon}\Big|\leq \Big|\frac{d(x,z)-d(y,z)}{d(x,y)}\Big|\leq 1.
	\end{align*}
	Hence, by~\cite[Theorem 3.1 and Theorem 4.1]{ks03}, for each $z\in X$, the energy measure of the distance function $d_z\colon x\mapsto d(x,z)$ satisfies
	\begin{align*}
		\int_X \varphi(x)d\Gamma_0^b(d_z,d_z)&=\lim_{\varepsilon\to 0}\frac{1}{b(\varepsilon)}\int_X \varphi(x)\int_{B(x,\varepsilon)}\Big(\frac{d(x,z)-d(y,z)}{\varepsilon}\Big)^2d\mu(y)d\mu(x)\\
		&\leq \lim_{\varepsilon\to 0}\int_X \varphi(x)\Big(\frac{\mu(B(x,\varepsilon))}{b(\varepsilon)}\Big)d\mu(x)\\
		&\leq \theta\int_X\varphi(x)d\mu(x).
	\end{align*}
	Since $\theta$ can be chosen to be arbitrarily close to 1, we conclude that 
	\begin{align*}
		\int_X \varphi(x)d\Gamma_0^b(d_z,d_z)\leq \int_X\varphi(x)d\mu(x).
	\end{align*}
	Consequently, $\Gamma_0^b(d_z,d_z)\leq \mu$ and so $d_0(x,x')\geq d_z(x)-d_z(x')$ for all $x,x',z\in X$.  The claim follows by selecting $z=x'$.
	
	\underline{Claim 2:} for each compact set $K\subset X$, there exists a constant $L$ such that $d_0\leq Ld$.
	
	Take $u\in \lip_0(X)$ for which $\Gamma_0^b(u,u)\leq \mu$. Our aim is to show that there exists a constant $L=L(K,c_\mu)$, depending only on $K$ and on the doubling constant $c_\mu$, such that for each $z\in K$, 
	$$\lip u(z)=\limsup_{\varepsilon\to 0}\sup_{x\in B(z,\varepsilon)}\frac{|u(x)-u(z)|}{\varepsilon}\leq L,$$
	from which we infer that $|u(x)-u(y)|\leq Ld(x,y)$ for all $x,y\in K$.
	
	Fix $\varepsilon>0$, $z\in K$ and $$L=L(\varepsilon,z)=\sup\Big\{\frac{|u(x)-u(y)|}{\varepsilon}:(x,y)\in B(z,\varepsilon)\times B(z,\varepsilon) \Big\}.$$
	%We need to show that for each $x\in X$, there exist $L=L_x>0$ and $r_x>0$ such that for each $\varepsilon<r_x$ and each $y\in B(x,\varepsilon)$, $|u(x)-u(y)|\leq Lr$. 
	%Fix $z\in X$ and $L=L(z)=\limsup_{\varepsilon\to 0}\sup_{y\in B(z,\varepsilon)}\frac{|u(z)-u(y)|}{\varepsilon}$.
	Then
	$$|u(x)-u(y)|\leq L\varepsilon\quad \text{ for all }x, y\in B(z,\varepsilon).$$
	Choose $(z_1,z_2)\in B(z,\varepsilon)\times B(z,\varepsilon)$ such that 
	$$|u(z_1)-u(z_2)|\geq \frac{1}{2}L\varepsilon>0.$$
	For $x\in B(z_1,\frac{1}{8}\varepsilon)$ and $y\in B(z_2,\frac{1}{8}\varepsilon)$,
	\begin{align*}
		|u(x)-u(y)|&\geq |u(z_1)-u(z_2)|-|u(z_1)-u(x)|-|u(y)-u(z_2)|\\
		&\geq \frac{1}{2}L\varepsilon-2L\varepsilon\frac{1}{8}=\frac{1}{4}L\varepsilon.
	\end{align*}
	
	Then the weak (2,2)-Poincar\'e inequality from~\cite[Theorem 4.2]{ks03} implies that for each $z\in K$ and $\varepsilon<R$, 
	\begin{align*}
		\int_{B(z,2\varepsilon)}\dashint_{B(z,2\varepsilon)}\Big(\frac{u(x)-u(y)}{\varepsilon}\Big)^2d\mu(y)d\mu(x)\leq C\mu(B(z,2\lambda\varepsilon)),
	\end{align*}
	for some constant $C$ depending on $K$ and $\lambda\geq 1$. Note that
	\begin{align*}
		\int_{B(z,2\varepsilon)}&\int_{B(z,2\varepsilon)}\Big(\frac{u(x)-u(y)}{\varepsilon}\Big)^2d\mu(y)d\mu(x)\\
		&\geq \int_{B(z_1,\varepsilon/8)}\int_{B(z_2,\varepsilon)}\Big(\frac{u(x)-u(y)}{\varepsilon}\Big)^2d\mu(y)d\mu(x)\\
		&\geq \frac{1}{16}L^2\mu(B(z_1,\varepsilon/8))\mu(B(z_2,\varepsilon/8)).
	\end{align*}
	Combining these two inequalities, we obtain
	\begin{align*}
		L^2\leq 16C\frac{\mu(B(z,2\varepsilon))\mu(B(z,\lambda \varepsilon))}{\mu(B(z_1,\varepsilon/8))\mu(B(z_2,\varepsilon/8))}\leq \hat{C},
	\end{align*}
	where $\hat{C}$ depends on the compact set $K$ and on the doubling constant of $\mu$. Thus we have shown that for each $z\in K$ and $\varepsilon>0$, there exists $L=L(K,c_\mu)$ such that
	$$|u(x)-u(y)|\leq L\varepsilon$$
	holds for all $x,y\in B(z,\varepsilon)$. In particular, this implies that $\lip u(z)\leq L$ for all $z\in K$ as desired.
	
	With these understood, property (C2) follows from the previous fact and~\cite[Theorem 4.2]{ks03}. 
	
	Finally, that $e_\varepsilon^{b,Y}(u)\rightharpoonup e^{b,Y}(u)$ follows from~\cite[Theorem 4.1]{ks03}.
\end{proof}

Next, we shall derive some useful consequences for approximating energy density to have (strong) property $\Cp$.

\bl\label{lemma:consequence of C}
Assume the approximating energy density $e_\varepsilon^c$ on a complete metric space $X$ has property $\Cp$. Then for each $u,v\in \mathbb{D}(\mathcal{E}_0)$ and $\varphi\in \lip_0(X)$, we have 
\begin{equation}\label{eq:condition 2 consequence}
	\lim_{\varepsilon\to 0}\int_X \varphi(x)	\frac{1}{c(x,\varepsilon)}\int_{B(x,\varepsilon)}\frac{(u(x)-u(x'))(v(x)-v(x'))}{\varepsilon^2}d\mu(x')d\mu(x)=\int_X \varphi(x)d\Gamma_0(u,v)(x).
\end{equation}  
\el
\begin{proof}
	We first show that for each $u\in \mathbb{D}(\mathcal{E}_0)$ and $\varphi\in \lip_0(X)$,
	$$\lim_{\varepsilon\to 0} E^c_\varepsilon(u)(\varphi)=\int_X \varphi d\Gamma_0(u,u),$$
	or more precisely, 
	\begin{equation}\label{eq:condition 2}
		\lim_{\varepsilon\to 0}\int_X \varphi(x)\frac{1}{c(x,\varepsilon)}	\int_{B(x,\varepsilon)}\frac{(u(x)-u(x'))^2}{\varepsilon^2}d\mu(x')d\mu(x)=\int_X \varphi(x)d\Gamma_0(u,u)(x).
	\end{equation}  
	Indeed, by~\eqref{eq:energy measure of Dirichlet form} and property (C1), we have
	\begin{align*}
		\int_X \varphi d\Gamma_0(u,u)&=\frac{1}{2}\big[2\mathcal{E}_0(u,\varphi u)-\mathcal{E}_0(u^2,\varphi) \big]\\
		&=\frac{1}{2}\Big[\lim_{\varepsilon\to 0}\int_{X}\frac{1}{c(x,\varepsilon)}\int_{B(x,\varepsilon)}\frac{2[u(x)-u(y)][\varphi(x)u(x)-\varphi(y)u(y)]}{\varepsilon^2}-\\
		&\qquad\qquad\qquad\qquad\qquad -\frac{(u^2(x)-u^2(y))(\varphi(x)-\varphi(y))}{\varepsilon^2}\mu(y)d\mu(x) \Big]\\
		&=\lim_{\varepsilon\to 0}\int_X \varphi(x)\frac{1}{c(x,\varepsilon)}	\int_{B(x,\varepsilon)}\frac{(u(x)-u(y))^2}{\varepsilon^2}d\mu(y)d\mu(x)+\\
		&\qquad\qquad\qquad\frac{1}{2}\lim_{\varepsilon\to 0}\int_X 	\frac{1}{c(x,\varepsilon)}\int_{B(x,\varepsilon)}(\varphi(y)-\varphi(x))\frac{(u(x)-u(y))^2}{\varepsilon^2}d\mu(y)d\mu(x)\\
		&=\lim_{\varepsilon\to 0} E^c_\varepsilon(u)(\varphi).
	\end{align*}

	Note that 
	\begin{align*}
		2(u(x)-u(x'))(v(x)-v(x'))&=\Big((u+v)(x)-(u+v)(x')\Big)^2-\\
		&\qquad\qquad\qquad (u(x)-u(x'))^2-(v(x)-v(x'))^2.
	\end{align*}
	The claim follows from~\eqref{eq:condition 2} together with the identity
	$$2\int_X \varphi d\Gamma_0(u,v)=\int_X \varphi d\Gamma_0(u+v,u+v)-\int_X \varphi d\Gamma_0(u,u)-\int_X \varphi d\Gamma_0(v,v).$$
\end{proof}

\bl\label{lemma:Absolute continuity of energy measure}
If the approximating energy functional $e_\varepsilon^{c,Y}$ between a complete metric space $X$ and an NPC space $Y$ has strong property $\Cp$, then for each $u\in KS^{1,2}(X,Y)$, $\mu_{u,Y}$ is absolutely continuous with respect to $\mu$.
\el 
\begin{proof}
	Fix a relatively compact domain $K$ in $X$. Suppose $\mu$ is uniformly doubling on balls with radius less than $r_0$. It suffices to show that there exists a constant $C$ (depending only on the local uniform doubling constant of $\mu$ on $K$ and on $K$) such that for each $B(x_0,r)\subset K$, $r\leq r_0$,
	\begin{align}\label{eq:aim}
	\int_{B(x_0,r)}d\mu_{u,Y}(x)\leq C\int_{B(x_0,r)}g_u^2d\mu(x).
	\end{align}
	Note that~\eqref{eq:aim} is isometrically invariant. Indeed, by the proof of Theorem 4.1 in~\cite{ks03}, the measure $\mu_{u,Y}$ remains un-changed if we isometrically embeds $Y$ to a Banach space (since the definition of energy is isometrically invariant). The same is true for the right-hand side of~\eqref{eq:aim}. Thus we may assume that $Y$ is isometrically contained in a Banach space with norm $\|\cdot\|$.
	
	Since
	$$\lim_{\varepsilon\to 0}\int_{B(x_0,r)}e_\varepsilon^{c,Y}(u)(x)d\mu(x)=\int_{B(x_0,r)}d\mu_{u,Y}(x).$$
	It suffices to give a uniform bound (independent of $\varepsilon$) on $\int_{B(x_0,r)}e_\varepsilon^{c,Y} (u)(x)d\mu(x)$ in terms of the right-hand side of~\eqref{eq:aim}.
	
	Let $g_u$ be the minimal 2-weak upper gradient of $u$. 
%	Then by the proof of Theorem 10.4.5 in~\cite{hkst12} (apply it with $X=B(x_0,r)$ and $p=2$ there), for $\varepsilon>0$ (and less than $r/2$) 
    Then by \eqref{eq:strong property C 3}, we have
	\begin{align*}
	\int_{B(x_0,r)}e_\varepsilon^{c,Y}(u)(x)d\mu(x)\leq C\int_{B(x_0,r)}g_u^2d\mu(x),
	\end{align*}
	where $C$ is independent of $\varepsilon$. 
	%depends only on the doubling constant of $\mu$ and the constant associated to the Poincar\'e inequality on $K$. 
	Consequently, 
	\begin{align*}
	\limsup_{\varepsilon\to 0}\int_{B(x_0,r)}e_\varepsilon^{c,Y}(u)(x)d\mu(x)\leq C\int_{B(x_0,r)}g_u^2d\mu(x),
	\end{align*}
	from which~\eqref{eq:aim} follows.
\end{proof}

In particular, if $X$ is uniformly locally doubling and supports a local uniform $(1,2)$-Poincar\'e inequality, then we would have 
$$d\mu_u(x)\leq C(\text{ap-Lip}u)^2(x)d\mu$$
for all $u\in KS^{1,2}(X,Y)$, where $C$ depends only on the associated data. Recall that for a map $h\colon X\to Z$, the approximate pointwise Lipschitz constant of $h$ is defined as
$$\text{ap-Lip}h(x):=\inf_{A}\limsup_{x'\to x,x'\in A\backslash\{x\}}\frac{d(h(x'),h(x))}{d(x',x)},$$
with the infimum taken over subsets $A\subset X$ having a Lebesgue point of density at $x$. In general, if $x$ is isolated, we let $\text{ap-Lip} u(x)=0$.

\bl\label{lemma:Poincare inequality for metric-valued mapping}
Assume the approximating energy functional $e_\varepsilon^{c,Y}$ between a complete metric space $X$ and an NPC space $Y$ has strong property $\Cp$. Then for each $u\in KS^{1,2}(X,Y)$, each compact set $K\subset (X,d_0)$ with $\lambda_KB\subset K$, the following Poincar\'e inequality holds:
$$\dashint_B\int_Bd^2(u(x),u(y))d\mu(x)d\mu(y)\leq C_K(\diam B)^2\int_{\lambda_K B}d\mu_{u,Y}(x).$$
\el
\begin{proof}
	Since $e_\varepsilon^c$ has property $\Cp$, by~\cite[Lemma 2.4 and Theorem 2.2]{kz12}, for each $v\in KS^{1,2}(X)$, we have the following (2,2)-Poincar\'e inequality,
	$$\int_B|v-v_B|^2d\mu(x)\leq C_K(\diam B)^2\int_{\lambda_K B}g_v^2d\mu.$$
	It is well-know that the previous inequality implies that for each $u\in KS^{1,2}(X,V)$, where $V=(V,\|\cdot\|_V)$ is a Banach space, we have 
	$$\dashint_B\int_B\|u(x)-u(y)\|_V^2d\mu(x)d\mu(y)\leq \hat{C}_K(\diam B)^2\int_{\lambda_K B}g_u^2d\mu.$$
	For this, see e.g.~\cite[Proof of Theorem 8.1.42]{hkst12} or~\cite[Proof of Theorem 3.6]{kst04}. On the other hand, \eqref{eq:strong property C 3} implies 
	%since $\mu$ is a locally doubling measure on $X$, Theorem~\ref{thm:KS is included in Newton} implies that 
	$$\int_{\lambda_K B}g_u^2d\mu\leq C_K\limsup_{\varepsilon\to 0}\int_{\lambda_K B}e_\varepsilon^{c,Y}(u)(x)d\mu.$$
	Consequently, for each $u\in KS^{1,2}(X,V)$, we have 
		$$\dashint_B\int_B\|u(x)-u(y)\|_V^2d\mu(x)d\mu(y)\leq C_K(\diam B)^2\int_{\lambda_K B}d\mu_{u,V}.$$
	Arguing as in the proof of Lemma~\ref{lemma:Absolute continuity of energy measure} (embedding $Y$ isometrically into some Banach space $V$; say $V=l^\infty(Y)$), we directly obtain that for each $u\in KS^{1,2}(X,Y)$
	$$\dashint_B\int_Bd^2(u(x),u(y))d\mu(x)d\mu(y)\leq C_K(\diam B)^2\int_{\lambda_K B}d\mu_{u,Y}(x).$$
\end{proof}

In our following proofs, the metric space $Y$ is typically fixed and  we will write $\mu_u$, instead of $\mu_{u,Y}$ for notational simplicity.

\subsection{Proof of Theorem~\ref{thm:main thm regularity} 1).}

We first show that the composition of the distance function with a harmonic mapping is subharmonic in the sense of~\cite{j97,bm95}. The proof is similar to that of Lemma 5 in~\cite{j97}; see also~\cite[Proof of Lemma 10.2]{ef01}.

\bp\label{prop:composition with distance subharmonic}
If $u\colon X\to Y$ is harmonic, then for each $y_0\in Y$ and each relatively compact open set $U\subset X$, the function $f_{y_0}:=d^2(u(\cdot),y_0)\colon X\to \R$ is weakly subharmonic on $U$, i.e., for each positive Lipschitz function $\lambda\colon X\to \R$ with $\supp(\lambda)\subset U$, 
\begin{equation*}
	-\mathcal{E}_0(\lambda,f_{y_0})=(\lambda,Af_{y_0})\geq 0.
\end{equation*}
Moreover, if $X$ has strong property C, then 
$$-\mathcal{E}_0(\lambda, f_{y_0})\geq 2\int_{X}\lambda d\mu_u(x).$$
\ep 

\begin{proof}
We write $v=d(u(\cdot),y_0)$. Note that both $v$ and $f_{y_0}(=v^2)$ are in $\mathbb{D}_{\loc}(\mathcal{E}_0)=N^{1,2}_{\loc}(X)$. Indeed, since the distance function is Lipschitz and $u\in KS^{1,2}_{\loc}(X,Y)$, we easily infer that $v\in KS^{1,2}_{\loc}(X)$. As, under our assumption on $X$, $KS^{1,2}_{\loc}(X)=N^{1,2}_{\loc}(X)$, we arrive at the desired observation. 

Let $\lambda \in \lip(X)$, $0\leq \lambda\leq 1$ be a Lipschitz function with $\supp(\lambda)\subset U$. For each $x\in X$, let $\gamma_x\colon [0,1]\to Y$ be the constant-speed geodesic in $Y$ from $\gamma_x(0)=u(x)$ to $\gamma_x(1)=y_0$. Define a map $u_\lambda\colon X\to Y$ by
\begin{align*}
	u_{\lambda}(x)=\gamma_x(\lambda(x)),\quad x\in X.
\end{align*}
Then $u_\lambda=u$ outside $U$. We next compare the energy of $u$ with that of $u_{\lambda}$.
 
Note that 
$$d(u_\lambda(x),y_0)\leq d(u(x),y_0)=v(x)$$
with $v\in L^2(X)$ since $u\in L^2(X,Y)$. This implies that $u_\lambda\in KS^{1,2}(X,Y)$.

Since $Y$ is NPC, triangle comparison (see \cite[Equation (2.1ii)]{ks93}) gives for $x,x'\in X$,
\begin{align*}
	d^2(u(x),u_\lambda(x'))&\leq (1-\lambda(x'))d^2(u(x),u(x'))+\lambda(x')d^2(u(x),y_0)\\
	&\quad -\lambda(x')(1-\lambda(x'))d^2(u(x'),y_0)
\end{align*} 
and 
\begin{align*}
d^2(u_\lambda(x),u_\lambda(x'))&\leq (1-\lambda(x))d^2(u(x),u(x'))+\lambda(x)d^2(u(x'),y_0)\\
&\quad -\lambda(x)(1-\lambda(x))d^2(u(x),y_0).
\end{align*} 
Inserting $d(u_\lambda(x'),y_0)=(1-\lambda(x'))d(u(x'),y_0)$, $d(u(\cdot),y_0)=v$ gives 

\begin{align*}
d^2&(u_\lambda(x),u_\lambda(x'))-d^2(u(x),u(x'))\\ 
&\leq -[\lambda(x)+\lambda(x')-\lambda(x)\lambda(x')]d^2(u(x),u(x'))\\
&\quad -(\lambda(x)-\lambda(x'))[(1-\lambda(x))f_y(x)-(1-\lambda(x'))f_y(x')].\numberthis\label{eq:key equation 1}
\end{align*} 
Since $|\lambda(x')-\lambda(x)|\leq L\cdot \varepsilon$ when $x'\in B(x,\varepsilon)$, we have
$$\lambda(x)+\lambda(x')-\lambda(x)\lambda(x')=2\lambda(x)-\lambda(x)^2+O(\varepsilon).$$
Now we have by~\eqref{eq:key equation 1}
\begin{align*}
	\frac{1}{c(x,\varepsilon)}\int_{B(x,\varepsilon)}&\frac{d^2(u_\lambda(x),u_\lambda(x'))-d^2(u(x),u(x'))}{\varepsilon^2}d\mu(x')\\ 
	&\leq -	\frac{1}{c(x,\varepsilon)}\int_{B(x,\varepsilon)}\frac{[2\lambda(x)-\lambda(x)^2+O(\varepsilon)]d^2(u(x),u(x'))}{\varepsilon^2}\\
	&\quad\quad -\frac{(\lambda(x)-\lambda(x'))[(1-\lambda(x))f_y(x)-(1-\lambda(x'))f_y(x')]}{\varepsilon^2}d\mu(x').
\end{align*}
Multiply by $\eta\in C_0(X)$, $0\leq \eta\leq 1$, on both side and integrate with respect to $x$, we get
\begin{align*}
\int\eta(x)\frac{1}{c(x,\varepsilon)}\int_{B(x,\varepsilon)}&\frac{d^2(u_\lambda(x),u_\lambda(x'))-d^2(u(x),u(x'))}{\varepsilon^2}d\mu(x')d\mu(x)\\ 
&\leq -	\int \frac{\eta(x)}{c(x,\varepsilon)}\int_{B(x,\varepsilon)}\frac{[2\lambda(x)-\lambda(x)^2+O(\varepsilon)]d^2(u(x),u(x'))}{\varepsilon^2}\\
&\qquad +\frac{(\lambda(x)-\lambda(x'))[(1-\lambda(x))f_y(x)-(1-\lambda(x'))f_y(x')]}{\varepsilon^2}d\mu(x')d\mu(x).
\end{align*}
Since $u$ is energy minimizing, the left-hand side of the above inequality is non-negative if we take $\limsup_{\varepsilon\to 0}$. Thus we infer that for $\varepsilon$ sufficiently small
\begin{align*}
	-\kappa(\varepsilon)&+\int \frac{\eta(x)}{c(x,\varepsilon)}\int_{B(x,\varepsilon)}\frac{(\lambda(x)-\lambda(x'))[(1-\lambda(x))f_y(x)-(1-\lambda(x'))f_y(x')]}{\varepsilon^2}d\mu(x')d\mu(x)\\
	&\leq -\int\frac{\eta(x)}{c(x,\varepsilon)}\int_{B(x,\varepsilon)}\frac{[2\lambda(x)-\lambda(x)^2+O(\varepsilon)]d^2(u(x),u(x'))}{\varepsilon^2}d\mu(x')d\mu(x),
\end{align*}
where $\kappa(\varepsilon)$ is a positive function in $\varepsilon$ that tends to 0 as $\varepsilon\to 0$. In the above inequality, replace $\lambda$ by $t\lambda$ with $t=\kappa(\varepsilon)^{1/2}$, divide by $t$ on both sides and then let $\varepsilon\to 0$, we obtain 
\begin{align*}
&\limsup_{\varepsilon\to 0}\int \frac{\eta(x)}{c(x,\varepsilon)}\int_{B(x,\varepsilon)}\frac{(\lambda(x)-\lambda(x'))(f_{y_0}(x)-f_{y_0}(x'))}{\varepsilon^2}d\mu(x')d\mu(x)\\
&\leq \limsup_{\varepsilon\to 0}-\int \frac{\eta(x)}{c(x,\varepsilon)}\int_{B(x,\varepsilon)}\frac{2\lambda(x)d^2(u(x),u(x'))}{\varepsilon^2}d\mu(x')d\mu(x)\leq 0.
\end{align*}
Together with Lemma~\ref{lemma:consequence of C}, we deduce that
\begin{align*}
\int_{X}\eta d\Gamma_0(\lambda,f_{y_0})\leq 0.
\end{align*}
In particular, we have $(\lambda,Af_{y_0})\geq 0$ for all positive Lipschitz function $\lambda$ with $\supp(\lambda)\subset U$. If $X$ has strong property $\Cp$, then we get from the previous inequality that
$$\int_X\eta d\Gamma_0(\lambda,f_{y_0})\leq -2\int_X \eta(x)\lambda(x)d\mu_u(x).$$
\end{proof}

We next point out that as in the setting of Riemannian manifolds, subharmonic functions are locally bounded.
\bl\label{lemma:boundedness of weakly harmonic maps}
Every positive weakly subharmonic function $v\colon X\to \R$ is locally bounded.
\el
\begin{proof}
By~\cite[Theorem 5.4]{bm95}, for every $p>0$, there exists a positive constant $c_p$ such that
\begin{align}\label{eq:local harnack}
	\sup_{B_0(x,r/2)}v\leq c_p\Big(\dashint_{B_0(x,r)}v^pd\mu(x)\Big)^{1/p}.
\end{align}
The claim follows by taking $p=2$.
%applying triangle inequality for the term on the right-hand side of~\eqref{eq:local harnack} (with $p=2$) together with the weak (2,2)-Poincar\'e inequality.
\end{proof}

Our proof of interior H\"older regularity follows the general approach of Lin~\cite[Proof of Theorem 3.1]{lin97}. For this, we need the following key covering type lemma, which generalizes~\cite[Lemma 3.5]{lin97} via a similar idea to the current setting.

From now on (till the end of this section), we assume that $Y$ is a locally doubling NPC space and the approximating energy functional $e_{\varepsilon}^{c,Y}$ between $X$ and $Y$ has property $\Cp$. Recall that a metric space $Y$ is doubling with constant $M$, $M\in \mathbb{N}$, if for each  ball $B(x,r)$, every $r/2$-separated subset of $B(x,r)$ has at most $M$ points. $Y$ is locally doubling if each compact subset $K$ of $Y$ is doubling with some constant $C_K$. 

\bl\label{lemma:covering}
Let $u\colon X\to Y$ be given as in Proposition~\ref{prop:composition with distance subharmonic}. We further assume that $\mu$ is doubling on $X$ with constant $c_d$, $X$ supports a weak $(2,2)$-Poincar\'e inequality with constant $C_P$ and $Y$ is $M$-doubling. Suppose $\kappa=\diam u(B(x_0,r))\in [\kappa_1,\kappa_0]$. There exists some $\varepsilon_0>0$, depending only on $\kappa_0,\kappa_1$ and $c_d$, $C_P$ and $M$ such that if $u(B(x_0,r))$ is covered by $k$ balls $B_1,\dots,B_k$ of radius $\varepsilon\leq \varepsilon_0$, then $u(B(x_0,r/2))$ can be covered by at most $k-1$ balls among $B_1,\dots,B_k$.	
\el

\begin{proof}
	For $i=1,2,\dots,k$, we take $x_i\in B(x_0,r)$ such that $B_i\subset B(p_i,2\varepsilon)$, where $p_i=u(x_i)$. Since $\varepsilon\leq \varepsilon_0\leq \kappa/16$, the balls $B(p_i,\kappa/8)$ covers $u(B(x_0,r))$. Since $u(B(x_0,r))$ has diameter $\kappa$, every $p_i$ belongs to a closed ball $\overline{B}'$ of radius $\kappa$ in $Y$. 
	
	Let $k'$ be the maximal number of points in the ball such that the distance is at least $\kappa/8$ apart. Since $Y$ is doubling, $k'\leq C_Y$ for some constant $C_Y$ depending only on the doubling constant $M$ of $Y$. Thus, we may assume that $B(p_i,\kappa/4)$, $i=1,2,\dots,k'$ covers $u(B(x_0,r))$. It follows that for at least one of those $p_i$, say for $p_1$, 
	\begin{align*}
		\mu\Big(u^{-1}\big(B(p_1,\kappa/4)\big)\cap B(x_0,r/2)\Big)\geq \frac{1}{k'}\mu(B(x_0,r/2))
	\end{align*}
	
	Consider the auxiliary function $f_{p_1}(x):=\frac{1}{\kappa^2}d^2(u(x),p_1)$. It is clear that
	\begin{align*}
		\tau:=\sup_{x\in B(x_0,r)}f_{p_1}(x)\leq \frac{1}{\kappa^2}\big(\diam u(B(x_0,r))\big)^2\leq 1.
	\end{align*}
	By the triangle inequality, and since $\diam u(B(x_0,r))=\kappa$, there exists some $\hat{x}\in B(x_0,r)$ with $f_{p_1}(\hat{x})\geq \kappa/2$, and so $\tau \geq \frac{1}{4}$. Since the function $h(x):=\tau-f_{p_1}(x)\geq 0$ on $B(x_0,r)$ and $h(x)\geq \frac{1}{8}$ on $u^{-1}\big(B(p_1,\kappa/4)\big)$. Since $f_{p_1}$ is weakly subharmonic, $h$ is weakly supharmonic and so by \cite[Lemma 6]{j97}, 
	\begin{align}\label{eq:inf estimate}
		\inf_{x\in B(x_0,r/2)}h(x)\geq c_1\dashint_{B(x_0,r)}h(x)d\mu(x)\geq \frac{c_1}{k'}\frac{\mu(B(x_0,r/2))}{\mu(B(x_0,r))}\geq c_0.
	\end{align}
	
	We next claim that for sufficiently small $\varepsilon$, we cannot have $u(B(x_0,\frac{r}{2}))\cap B(p_i,2\varepsilon)\neq \emptyset$ for all $i=1,\dots,k$. Indeed, let $\tilde{x}\in B(x_0,r)$ be such that $\tau=\frac{1}{\kappa^2}d^2(u(\tilde{x}),p_1)$. since the balls $B(p_i,2\varepsilon)$ cover $u(B(x_0,r))$, we find some $p_i$ with $d(u(\tilde{x}),p_i)\leq 2\varepsilon$. Thus, if $u(B(x_0,\frac{r}{2}))\cap B(p_i,2\varepsilon)\neq \emptyset$ for all $i=1,\dots,k$, we would have $d(u(x_1),u(\tilde{x}))\leq 4\varepsilon$ for some $x_1\in B(x_0,r/2)$ and so 
	\begin{align*}
		\inf_{x\in B(x_0,r/2)}h(x)&\leq h(x_1)=\tau-\frac{1}{\kappa^2}d^2(u(x_1),p_1)\\
		&\leq \tau-\frac{1}{\kappa^2}(d^2(u(\tilde{x}),p_1)+d^2(u(\tilde{x}),u(x_1))-2d(u(\tilde{x}),p_1)d(u(\tilde{x}),u(x_1)))\\
		&\leq\frac{16\varepsilon\sqrt{\tau}}{\kappa},
	\end{align*}  
	which contradicts~\eqref{eq:inf estimate} if $\varepsilon_0$ is sufficiently small, quantitatively.
\end{proof}

With the help of Lemma~\ref{lemma:covering}, the proof of interior H\"older regularity can be deduced by a similar argument as in Lin~\cite[Proof of Theorem 3.1]{lin97}.

\bp\label{prop:local regularity under PI}
Let $u\colon X\to Y$ be given as in Proposition~\ref{prop:composition with distance subharmonic}, then $u$ is locally H\"older continuous. 
\ep 
\begin{proof}
	Since the issue is local, we may assume that $\mu$ is a global doubling measure on $X$ with doubling constant $c_d$ and $X$ supports a weak $(2,2)$-Poincar\'e inequality with constant $C_P$ and $Y$ is $M$-doubling. In below, we refer $c_d$, $C_P$ and $M$ as the data of $X$ and $Y$. Under these assumptions, we first show that there exists $\delta_0$, depending only on the data of $X$ and $Y$ such that if $B(x_0,2)\subset\subset \Omega$, then $\diam u(B(x_0,1))=2$ implies that $\diam u(B(x_0,\delta))\leq 1$. For notational simplicity, we write $B_t=B(x_0,t)$ for $t>0$. 
	
	Let $\varepsilon_0$ be given as in Lemma~\ref{lemma:covering}. Since $\diam u(B_1)=2$, we can cover $u(B_1)$ by $k$ balls of radius $\varepsilon_0$ in $N$, where $k$ depends only on the data of $X$ and $Y$. By Lemma~\ref{lemma:covering}, $u(B_{2^{-1}})$ can be covered by at most $k-1$ balls of radius $\varepsilon_0$. If $\diam u(B_{2^{-1}})>1$, we may repeat the above arguments with  $2^{-1}$ in place of 1 and $k-1$ in place of $k$ to conclude that $u(B_{2^{-2}})$ can be covered by at most $k-2$ balls. It follows that there is some $k_0\leq k$ such that $\diam u(B_{2^{-k_0}})\leq 1$. In particular, we may set $\delta=2^{-k_0}$ for the claim.
	
	Consider now $B(x_0,r)$ with $B(x_0,2r)\subset\subset \Omega$. Let $\kappa=\diam u(B(x_0,r))$. Denote by $\hat{X}=(X,r^{-1}d,\mu)$ and $Y_\kappa=(Y,2\kappa^{-1} d)$. Then for each $v\in \lip_0(X)$, 
	\begin{align*}
		\hat{E}_\varepsilon(v)&=\int_{X}\frac{1}{c(x,r\varepsilon)}\int_{\hat{B}(x,\varepsilon)}\frac{|v(x)-v(y)|^2}{\varepsilon^2}d\mu(y)d\mu(x)\\
		&=\int_X \frac{1}{c(x,r\varepsilon)}\int_{B(x,r\varepsilon)}\frac{|v(x)-v(y)|^2}{\varepsilon^2}d\mu(y)d\mu(x)\\
		&=r^{2}E_{r\varepsilon}(v).
	\end{align*}
	Thus, $\lim_{\varepsilon\to 0}\hat{E}_\varepsilon(u)=r^{2}E_0(u)$ exists and induces a regular strongly local Dirichlet form $\hat{\mathcal{E}}_0=r^{2}\mathcal{E}_0$. Let $\hat{\Gamma}_0$ be the energy measure associated to the Dirichlet form $\hat{\mathcal{E}}_0$. Then $d\hat{\Gamma}_0=r^{2}d\Gamma_0$. 
	
	We next show that $\hat{X}$ and $Y_{\kappa}$ have the same data as that of $X$ and $Y$. Indeed, for each $\hat{B}(x,s)\subset \hat{X}$,
	\begin{align*}
		\mu(\hat{B}(x,2s))=\mu(B(x,2sr))\leq c_d\mu(B(x,sr))=c_d\mu(\hat{B}(x,s))
	\end{align*}
	and 
	\begin{align*}
		\int_{\hat{B}(x,s)}|v(y)-v_{\hat{B}(x,s)}|^2d\mu(y)&=\int_{B(x,sr)}|v-v_{B(x,sr)}|^2d\mu\\
		&\leq C_P (\diam B(x,sr))^2\int_{\lambda B(x,sr)}d\Gamma_0(v,v)\\
		&=C_p(\diam \hat{B}(x,s))^2\int_{\lambda \hat{B}(x,s)}d\hat{\Gamma}_0(v,v). 
	\end{align*}
	That $Y_\kappa$ is $M$-doubling is clear.
	
	We now regard our harmonic mapping $u\colon X\to Y$ as a mapping $\hat{u}\colon \hat{X}\to (Y,2\kappa^{-1} d)$. Then a simple computation implies 
	\begin{align*}
		\hat{E}(u_r)=\inf_{v\in KS^{1,2}_\phi(\Omega,Y)}\hat{E}(v),
	\end{align*}
	which means that $\hat{u}\colon \hat{X}\to Y_\kappa$ is harmonic. Moreover, 
	$$\diam_{Y_\kappa} \hat{u}(\hat{B}(x_0,1))=2k^{-1}\diam u(B(x_0,r))=2.$$
	Thus, we may apply Lemma~\ref{lemma:covering} to find some $\delta$, depending only on the data of $\hat{X}$ and $Y_{\kappa}$ (and hence on the data of $X$ and $Y$), such that $\diam_{Y_\kappa}\hat{u}(\hat{B}(x_0,\delta))\leq 1$. That is, $2k^{-1}\diam u(B(x_0,\delta r))\leq 1$, or equivalently, $\diam u(B(x_0,\delta r))\leq \frac{1}{2}\diam u(B(x_0,r))$. A standard iteration then gives the desired H\"older continuity.
\end{proof}

\begin{proof}[Proof of Theorem~\ref{thm:main thm regularity} 1)]
	This follows from Proposition~\ref{prop:local regularity under PI}.	
\end{proof}

\subsection{Proof of Theorem~\ref{thm:main thm regularity} 2)}

We will now follow the approach of Jost~\cite{j97} to prove the interior regularity for harmonic mappings from metric spaces with property $\Cp$ to NPC spaces. Throughout this section, we consider everything in the metric space $(X,d_0)$. For instance, when we say $B(x_0,R)\subset X$ is a ball, we mean that $B(x_,R)$ is a ball in $(X,d_0)$.

If $B(x_0,R)\subset X$ is a ball with fixed center $x_0$, we write
\begin{align*}
	v_{+,R}:=\sup_{x\in B(x_0,R)}v(x),\ v_{-,R}:=\inf_{x\in B(x_0,R)}v(x) \text{ and } v_R:=\dashint_{B(x_0,R)}v(x)d\mu(x).
\end{align*}

The following lemma is a variant of~\cite[Lemma 8]{j97}.

\bl\label{lemma:upper estimate for regularity}
Suppose the approximating energy functional $e_\varepsilon^{c,Y}$ between a complete metric space $X$ and an NPC space $Y$ has strong property $\Cp$. Let $u\colon X\to Y$ be a harmonic mapping. Fix $y_0\in Y$ and $B(x_0,4R)\subset\subset X$. Then 
\begin{align*}
\frac{R^2}{\mu(x_0,R)}\int_{B(x_0,R)}d\mu_u(x)\leq c(v_{+,4R}-v_{+,R}),
\end{align*}
where $v(x)=f_{y_0}(x)=d^2(u(x),y_0)$.
\el
\begin{proof}
Note that by Proposition~\ref{prop:composition with distance subharmonic} we have $2\int \lambda(x)d\mu_u(x)\leq -\mathcal{E}_0(\lambda,v)$ for each positive Lipschitz function $\lambda$ with compact support.

Fix $B(x_0,4R)\subset\subset X$. Let $G^R(x_0,x)$ be the singular Green function on $B(x_0,R)$ relative to $B(x_0,2R)$, that is, $G^R(x_0,\cdot)\in \mathbb{D}(\mathcal{E}_0,B(x_0,2R))$ (closure of $\mathbb{D}(\mathcal{E}_0)\cap C_0(B(x_0,2R))$) and it is the unique solution of
\begin{align*}
(-\int_{B(x_0,2R)}\varphi(x)AG^R(x_0,x)d\mu(x)=)		\mathcal{E}_0(\varphi,G^R(x_0,\cdot))=\dashint_{B(x_0,R)}\varphi(x)d\mu(x)
\end{align*}
for all $\varphi\in \mathbb{D}(\mathcal{E}_0)$ with $\supp(\varphi)\subset B(x_0,2R)$; see~\cite[Section 6]{bm95} for the existence and basic properties of this Green function. Set
$$w^R(x):=\frac{\mu(B(x_0,R))}{R^2}G^R(x_0,x)\in \mathbb{D}(\mathcal{E}_0,B(x_0,2R)).$$
Then 
\begin{align}\label{eq:6.1}
	\mathcal{E}_0(\varphi,w^R)=\frac{1}{R^2}\int_{B(x_0,R)}\varphi(x)d\mu(x)
\end{align}
for all $\varphi\in \mathbb{D}(\mathcal{E}_0)$ with $\supp(\varphi)\subset B(x_0,2R)$. Furthermore, by the estimates for $G^R$~\cite[Theorem 6.1]{bm95}, we have
\begin{align*}
	&0\leq w^R\leq \gamma_1\quad \text{ in } B(x_0,2R)\numberthis\label{eq:6.2} \\
	&w^R\geq \gamma_2>0\quad \text{ in }B(x_0,R)\numberthis\label{eq:6.3}
\end{align*}
for some structural constants $\gamma_1,\gamma_2$ that does not depend on $R$.

Set $z(x):=v(x)-v_{+,4R}$. Then we have
\begin{align*}
	2\int_{B(x_0,2R)}\big(w^R(x)\big)^2d\mu_u(x)&\leq -\mathcal{E}_0(\big(w^R\big)^2,z)\\
	&= -2\mathcal{E}_0(w^R,w^Rz)+2\int zd\Gamma_0(w^R,w^R)\\
	&\leq -2\mathcal{E}_0(w^R,w^Rz) \quad\text{since }z\leq 0.
\end{align*} 
From~\eqref{eq:6.1},~\eqref{eq:6.2},~\eqref{eq:6.3} and~\cite[Corollary 1]{j97}, we obtain
\begin{align*}
	\int_{B(x_0,R)}d\mu_u(x)&\leq c_1\int_{B(x_0,2R)}\big(w^R(x)\big)^2d\mu_u(x)\\
	&\leq -\frac{c_1}{R^2}\int_{B(x_0,R)}w^R(x)z(x)d\mu(x)\\
	&\leq \frac{c_2}{R^2}\int_{B(x_0,R)}v_{+,4R}-v(x)d\mu(x)\\
	&= \frac{c_2\mu(B(x_0,R))}{R^2}(v_{+,4R}-v_R)\\
	&\leq \frac{c_3\mu(B(x_0,R))}{R^2}(v_{+,4R}-v_{+,R}).
\end{align*}

\end{proof}

\begin{proof}[Proof of Theorem~\ref{thm:main thm regularity} 2)]
	With Lemma~\ref{lemma:Poincare inequality for metric-valued mapping} and Lemma~\ref{lemma:upper estimate for regularity} at hand, the proof of~\cite[Theorem, Section 6]{j97} works with minor changes in our setting.
	
	For each $\rho>0$, we set
	$$\omega(\rho):=\sup_{x\in B(x_0,\rho)}d^2(u(x),p)=v_{p,+,\rho}.$$
	Our aim is to show that for each $p$ in the convex hull of $u(B(x_0,\delta R))$, where $\delta$ is a fixed constant, and for all $\rho<\frac{R}{2}$ sufficiently small,
	$$\omega(\rho)\leq c\Big(\frac{\rho}{R}\Big)^\alpha \omega(R)$$
	for some constant $c>0$ and some $\alpha\in (0,1)$. 
	
	Before turn to the proof of the claim, we observe the claim implies the H\"older continuity via a standard argument as follows. Take $p=\bar{u}_\rho$, the mean value (or the center of mass) of $u$ on $B(x_0,\rho)$, then 
	$$\omega(\rho)^{1/2}\leq \text{osc}_{B(x_0,\rho)}u\leq 2\omega(\rho)^{1/2},$$
	from which the H\"older continuity follows.
	
	Since the issue is isometrically invariant, we may assume that $Y$ is (isometrically) contained in some Banach space $V$ with norm $\|\cdot\|$. Set $\bar{u}_R$ be the mean value of $u$ on $B(x_0,R)$ and 
	$v_p(x):=\|u(x)-p\|^2$. We shall apply~\cite[Lemma 7]{j97} to the function $v_{\bar{u}_{\frac{R}{4}}}$. Choose $\varepsilon=\frac{1}{8}$ and $R'\in [\varepsilon^mR,\frac{R}{4}]$. 
	%with $\varepsilon^mR\leq R'\leq \frac{R}{4}$. 
	Since $X$ has strong property C, we have by Lemma~\ref{lemma:Poincare inequality for metric-valued mapping} and Lemma~\ref{lemma:upper estimate for regularity} that 
	\begin{align*}
		v_{R'}&=\dashint_{B(x_0,R')}\|u(x)-\bar{u}_{\frac{R}{4}}\|^2d\mu(x)\leq C_0\dashint_{B(x_0,R/4)}\|u(x)-\bar{u}_{\frac{R}{4}}\|^2d\mu(x)\\
		&\leq \frac{C_1R^2}{\mu(B(x_0,R))}\int_{B(x_0,\lambda R)}d\mu_u(x)\leq C_2(v_{p,+,\lambda R}-v_{p,+,\lambda R/4}).
	\end{align*}
	Combining this estimate with~\cite[Lemma 7]{j97}, we get for each $p$ in the convex hull of $u(B(x_0,\delta R))$, $\delta=\varepsilon^m$, we have
	\begin{align*}
		\sup_{x\in B(x_0,\delta R)}\|u(x)-p\|^2&\leq 4\sup_{x\in B(x_0,\delta R)}\|u(x)-\bar{u}_{\frac{R}{4}}\|^2\\
		&\leq 4\varepsilon^2\sup_{x\in B(x_0,R)}\|u(x)-\bar{u}_{\frac{R}{4}}\|^2+C_3(v_{p,+,\lambda R}-v_{p,+,\lambda R/4})\\
		&\leq 16\varepsilon^2 \sup_{x\in B(x_0,R)}\|u(x)-\bar{u}_{\frac{R}{4}}\|^2+C_3(v_{p,+,\lambda R}-v_{p,+,\delta R}).
	\end{align*} 
	Since 
	$$\sup_{x\in B(x_0, R)}\|u(x)-p\|^2\leq \sup_{x\in B(x_0,\lambda R)}\|u(x)-p\|^2=\omega(\lambda R),$$
	we have
	$$(1+C_3)\omega(\delta R)\leq \Big(\frac{1}{64}+C_3\Big)\omega(\lambda R).$$
	A simple iteration then gives our desired estimate
	$$\omega(\rho)\leq c\Big(\frac{\rho}{R}\Big)^\alpha \omega(R)$$
	for some constant $c>0$ and some $\alpha\in (0,1)$. 
\end{proof}

\section{A Liouille type theorem for harmonic mappings}\label{sec:Liouville}

%\subsection{A review of Sobolev functions on quasiopen subsets of metric spaces}

\subsection{Composition with distance function}

We first show that the composition of the distance function with a harmonic mapping $u$ is ``almost'' weakly subharmonic. The proof is similar to that used in Proposition~\ref{prop:composition with distance subharmonic}, relying on the idea of Jost~\cite{j97}.

\bl\label{lemma:for composition with distance}
If $u\colon X\to Y$ is harmonic, then for each $y_0\in Y$, the function $v=d(u(\cdot),y_0)\colon X\to \R$ satisfies that for each positive $\lambda\in N^{1,2}(X)\cap L^\infty(X)$ with compact support, 
\begin{equation}\label{eq:for composition with distance}
-\mathcal{E}_0(\lambda v,v)\geq 0.
\end{equation}
\el

\begin{proof}
	Let $\lambda \in \lip_0(X)$, $0\leq \lambda\leq 1$ be a Lipschitz function with $\supp(\lambda)\subset U\subset\subset X$. By standard approximation, it suffices to prove \eqref{eq:for composition with distance} holds for all such $\lambda$. For each $x\in X$, let $\gamma_x\colon [0,1]\to Y$ be the constant-speed geodesic in $Y$ from $\gamma_x(0)=u(x)$ to $\gamma_x(1)=y_0$. Define a map $u_\lambda\colon X\to Y$ by
	\begin{align*}
	u_{\lambda}(x)=\gamma_x(\lambda(x)),\quad x\in X.
	\end{align*}
	Then $u_\lambda=u$ outside $U$. We next compare the energy of $u$ with that of $u_{\lambda}$.
	
	Note that 
	$$d(u_\lambda(x),y_0)\leq d(u(x),y_0)=v(x)$$
	with $v\in L^2(X)$ since $u\in L^2(X,Y)$. This implies that $u_\lambda\in KS^{1,2}(X,Y)$.
	
	Since $Y$ is NPC, triangle comparison (see \cite[Equation (2.1ii)]{ks93}) gives for $x,x'\in X$,
	\begin{align*}
	d^2(u(x),u_\lambda(x'))&\leq (1-\lambda(x))d^2(u(x),u(x'))+\lambda(x')d^2(u(x),y_0)\\
	&\quad -\lambda(x')(1-\lambda(x'))d^2(u(x'),y_0)
	\end{align*} 
	and 
	\begin{align*}
	d^2(u_\lambda(x),u_\lambda(x'))&\leq (1-\lambda(x))d^2(u(x),u(x'))+\lambda(x)d^2(u(x'),y_0)\\
	&\quad -\lambda(x)(1-\lambda(x))d^2(u(x),y_0).
	\end{align*} 
	Inserting $d(u_\lambda(x'),y_0)=(1-\lambda(x'))d(u(x'),y_0)$, $d(u(\cdot),y_0)=v$ gives 
	\begin{align*}
	d^2&(u_\lambda(x),u_\lambda(x'))-d^2(u(x),u(x'))\\ 
	&\leq -[\lambda(x)+\lambda(x')-\lambda(x)\lambda(x')][d^2(u(x),u(x'))-(v(x)-v(x'))^2]\\
	&\quad -2(v(x)-v(x'))[\lambda(x)v(x)-\lambda(x')v(x')]+[\lambda(x)v(x)-\lambda(x')v(x')]^2.\numberthis\label{eq:key equation}
	\end{align*} 
	Since $|\lambda(x')-\lambda(x)|\leq L\cdot \varepsilon$ when $x'\in B(x,\varepsilon)$, we have
	$$\lambda(x)+\lambda(x')-\lambda(x)\lambda(x')=2\lambda(x)-\lambda(x)^2+O(\varepsilon).$$
	Now we have by~\eqref{eq:key equation}
	\begin{align*}
	\frac{1}{c(x,\varepsilon)}\int_{B(x,\varepsilon)}&\frac{d^2(u_\lambda(x),u_\lambda(x'))-d^2(u(x),u(x'))}{\varepsilon^2}d\mu(x')\\ 
	&\leq -	\frac{1}{c(x,\varepsilon)}\int_{B(x,\varepsilon)}\frac{[2\lambda(x)-\lambda(x)^2+O(\varepsilon)][d^2(u(x),u(x'))-(v(x)-v(x'))^2]}{\varepsilon^2}\\
	&\quad +\frac{2(v(x)-v(x'))[\lambda(x)v(x)-\lambda(x')v(x')]+[\lambda(x)v(x)-\lambda(x')v(x')]^2}{\varepsilon^2}d\mu(x').
	\end{align*}
	Multiply by $\eta$ on both side and integrate with respect to $x$, we get
	\begin{align*}
	\int\frac{\eta(x)}{c(x,\varepsilon)}&\int_{B(x,\varepsilon)}\frac{d^2(u_\lambda(x),u_\lambda(x'))-d^2(u(x),u(x'))}{\varepsilon^2}d\mu(x')d\mu(x)\\ 
	&\leq -	\int \frac{\eta(x)}{c(x,\varepsilon)}\int_{B(x,\varepsilon)}\frac{[2\lambda(x)-\lambda(x)^2+O(\varepsilon)][d^2(u(x),u(x'))-(v(x)-v(x'))^2]}{\varepsilon^2}\\
	&\quad +\frac{2(v(x)-v(x'))[\lambda(x)v(x)-\lambda(x')v(x')]+[\lambda(x)v(x)-\lambda(x')v(x')]^2}{\varepsilon^2}d\mu(x')d\mu(x).
	\end{align*}
	Since $u$ is energy minimizing, the left-hand side of the above inequality is non-negative if we take $\limsup_{\varepsilon\to 0}$. Thus, we infer that 
	\begin{align*}
	-\kappa(\varepsilon)&+\int \frac{\eta(x)}{c(x,\varepsilon)}\int_{B(x,\varepsilon)}\frac{[2\lambda(x)-\lambda(x)^2][d^2(u(x),u(x'))-(v(x)-v(x'))^2]}{\varepsilon^2}d\mu(x')d\mu(x)\\
	&\leq -\int \frac{\eta(x)}{c(x,\varepsilon)}\int_{B(x,\varepsilon)}\frac{2(v(x)-v(x'))[\lambda(x)v(x)-\lambda(x')v(x')]}{\varepsilon^2}d\mu(x')d\mu(x),
	\end{align*}
	where $\kappa(\varepsilon)$ is a positive function in $\varepsilon$ that tends to $0$ as $\varepsilon\to 0$. In the above inequality, replace $\lambda$ by $t\lambda$ with $t=\kappa(\varepsilon)^{1/2}$, divide by $t$ and then let $\varepsilon \to 0$, we obtain 
	\begin{align*}
	&\limsup_{\varepsilon\to 0}\int \frac{\eta(x)}{c(x,\varepsilon)}\int_{B(x,\varepsilon)}\frac{\lambda(x)[d^2(u(x),u(x'))-(v(x)-v(x'))^2]}{\varepsilon^2}d\mu(x')d\mu(x)\\
	&\leq \limsup_{\varepsilon\to 0}-\int \frac{\eta(x)}{c(x,\varepsilon)}\int_{B(x,\varepsilon)}\frac{(v(x)-v(x'))[\lambda(x)v(x)-\lambda(x')v(x')]}{\varepsilon^2}d\mu(x')d\mu(x)\\
	&=-\int_X \eta(x)d\Gamma_0(\lambda v,v)(x),
	\end{align*}
	where in the last step we have applied Lemma~\ref{lemma:consequence of C}. Since $|v(x)-v(x')|\leq d(u(x),u(x'))$, the first term in the above inequality is non-negative, this proves the claim.
\end{proof}

We would like to point out the following Liouville theorem, which follows directly from Lemma \ref{lemma:for composition with distance}.

\bt\label{thm:Liouville subharmonic function} 
Let $X$ be a complete non-compact metric measure space and $Y$ an NPC space. Suppose the approximating energy density $e_\varepsilon^{c,Y}$ between $X$ and $Y$ has strong property $\Cp$. If $u\in KS^{1,2}_{\loc}(X,Y)\cap L^2(X,Y)$ is a harmonic mapping, then $u$ is constant.
\et 

\begin{proof}[Proof of Theorem~\ref{thm:Liouville subharmonic function}]
	Fix a base point $x_0\in X$ and define $\rho\colon X\to \R$ as 
	\begin{align*}
	\rho(x)=\max\Big\{0,\min\big\{1,2-\frac{1}{R}d(x,x_0) \big\} \Big\}.
	\end{align*}
	Then $\rho$ is $\frac{1}{R}$-Lipschitz and $\rho=0$ on $X\backslash B(x_0,2R)$ and $\rho=1$ on $B(x_0,R)$. For notational simplicity, we write $B_{r}=B(x_0,r)$. Since $v=d(u(\cdot),u(x_0))$ satisfies \eqref{eq:for composition with distance}, 
	\begin{align}\label{eq:eq 1}
	0\geq \mathcal{E}_0(\rho^2v, v)=\frac{1}{2}\mathcal{E}_0(\rho^2,v^2)+\int_X \rho^2d\Gamma_0(v,v).
	\end{align}
	By the Leibniz rule and Cauchy-Schwartz inequality for Dirichlet forms, we have
	\begin{align*}
	\mathcal{E}_0(\rho^2,v^2)^2&=\Big(\int_X d\Gamma_0(\rho^2,v^2) \Big)^2=16\Big(\int_X \rho vd\Gamma_0(\rho,v) \Big)^2\\
	&\leq \int_Xv^2d\Gamma_0(\rho,\rho)\cdot \int_X \rho^2d\Gamma_0(v,v). 
	\end{align*}
	Consequently, we have
	\begin{align*}
	\mathcal{E}_0(\rho^2,v^2)\geq -4\Big(\int_Xv^2d\Gamma_0(\rho,\rho)\Big)^{1/2}\cdot \Big(\int_X \rho^2d\Gamma_0(v,v)\Big)^{1/2}.
	\end{align*}
	Substituting the above estimate in~\eqref{eq:eq 1}, we obtain
	\begin{align*}\label{eq:2}
	0&\geq \int_{B_{2R}\backslash B_R}\rho^2 d\Gamma_0(v,v)-2\Big(\int_{B_{2R}\backslash B_R}v^2d\Gamma_0(\rho,\rho)\Big)^{\frac{1}{2}}\cdot \Big(\int_{B_{2R}\backslash B_R} \rho^2d\Gamma_0(v,v)\Big)^{\frac{1}{2}}\\
	&\qquad\qquad+\int_{B_R}d\Gamma_0(v,v),
	\end{align*}
	which is a polynomial, $P(\Psi)=\Psi^2-2b\Psi+c$ with $\Psi=\Big(\int_{B_2R\backslash B_R} \rho^2d\Gamma_0(v,v)\Big)^{\frac{1}{2}}$. Since it has non-positive value, it must hold $b^2\geq c$ and so 
	\begin{align*}
	\int_{B_R}d\Gamma_0(v,v)\leq \int_{B_{2R}\backslash B_R}v^2d\Gamma_0(\rho,\rho)\lesssim \frac{1}{R^2}\int_{B_{2R}}v^2d\mu.
	\end{align*} 
	In particular, 
	\begin{align*}
	\int_{B_R}d\Gamma_0(v,v)\lesssim \frac{1}{R^2}\int_{X}v^2d\mu.
	\end{align*}
	Since $v\in L^2(X)$, sending $R$ to infinity, we conclude that
	\begin{align*}
	\int_X d\Gamma_0(v,v)=\mathcal{E}_0(v,v)=0.
	\end{align*}	
	It follows from the Poincar\'e inequality that $v$ is constant and hence $u\equiv u(x_0)$  on $X$.
	
\end{proof}

\subsection{Proof of Theorem \ref{thm:Liouville}}

\bp\label{prop:composition with distance subharmonic 2}
If $u\colon X\to Y$ is harmonic, then for each $y_0\in Y$, the function $v=d(u(\cdot),y_0)\colon X\to \R$ is weakly subharmonic, i.e. for each positive Lipschitz function $\lambda\colon X\to \R$ with compact support $U\subset\subset X$, 
\begin{equation*}
-\mathcal{E}_0(\lambda,v)\geq 0.
\end{equation*}
\ep 

\begin{proof}
The passage from \eqref{eq:for composition with distance} to subharmonicity is well-known and relies on the fine topology and potential theory; see for instance in \cite{f05,ef01}. For the convenience of the readers, we present the proof here.

Let $\psi\in N^{1,2}(X)\cap L^\infty(X)$ be nonnegative function with $\supp(\psi)\subset U\subset \subset X$ and let $T_n\colon \R\to \R$ be the nearest point project of $\R$ onto $[\frac{1}{n},n]$. Then the cut-off functions $t_n:=T_n\circ d(\cdot,y_0)$ and $\frac{1}{t_n}$ are $n$-Lipschitz. Thus the truncated functions $v_n:=t_n\circ u$ and $\frac{1}{v_n}$ are in the Sobolev space $N_{loc}^{1,2}(X)$. Consider the sets 
$$U_n:=\Big\{x\in U:\frac{1}{n}<v(x)<n\Big\},\quad n\in \N.$$
For each $n\in \N$, $U_n$ is quasiopen as $v$ is quasicontinuous (see Section \ref{subsec:Sobolev spaces}). Thus $\{U_n\}_{n\in \N}$ forms a  covering of $U$ with quasiopen sets. 

By Lemma \ref{lemma:Kilpelainen-Maly}, we may find a sequence of nonnegative functions $\{\psi_j\}_{j\in N}$ converging to $\psi$ in $N^{1,2}_0(U)$ such that each $\psi_j$ is a finite sum of functions in $N^{1,2}_0(U_\alpha)$. Without loss of generality, we may assume each $\psi_j$ lies in $N^{1,2}_0(U_j)$.
%That is, $\psi_j=\sum_{k=1}^{N_j}\psi_{j,k}$ for $\psi_{j,k}\in N^{1,2}_0(U_{k})$.
Then $\lambda_n:=\frac{\psi_{n}}{v_n}\in N^{1,2}_0(U)\cap L^\infty(U)$ and $\psi_n=\lambda_n v$ in $U_n$. Moreover, by \eqref{eq:for composition with distance}, we have
$$-\mathcal{E}_0(\psi,v)=-\lim_{n\to\infty}\mathcal{E}_0(\psi_n,v)=-\lim_{n\to\infty}\mathcal{E}_0(\lambda_n v,v)\geq 0.$$
The proof is complete.
%This implies that $v$ is finely subharmonic off a polar subset of $U_n$ and hence q.e. in $U=\bigcup_nU_n$ as any countable union of polar sets is again polar.
\end{proof}

\begin{proof}[Proof of Theorem~\ref{thm:Liouville}]
	By Proposition~\ref{prop:composition with distance subharmonic 2}, the function $v(x)=d(u(x),u(x_0))$, $x_0\in X$, is weakly subharmonic. By~\cite[Theorem 1]{s94}, $v$ is constant and so $u$ is constant as well.
\end{proof}

\section{Harmonic mapping flow}\label{sec:harmonic mapping flow}
%In this section, we say that a metric space $X$ has strong property B if $X=(X,d,\mu)$ is locally uniformly doubling, supports a weak (1,2)-Poincar\'e inequality, and satisfies the properties (B2) and (B3) in the definition of property B. 
Throughout this section, we assume that $X$ is $\RCD(K,N)$ and $Y$ is NPC. Recall that the distance $D$ on $L^2(X,Y)$ is defined as
\begin{align*}
	D^2(u,v)=\int_X d^2(u(x),v(x))d\mu(x).
\end{align*}

The proof of the following elementary lemma can be found in~\cite[Corollary 4.1.1]{j97book}.
\bl\label{lemma:L2 is NPC}
	$L^2(X,Y)$ is NPC.
\el
%\begin{proof}
%We first show that	$L^2(X,Y)$ is a length space. Fix $u_0,u_1\in L^2(X,Y)$. For each $x\in X$, let $\alpha_x$ be the (unique) geodesic connecting $u_0(x)$ and $u_1(x)$. Then we define a curve $\gamma_t$ on $L^2(X,Y)$ as follows: for each $t$, $\gamma_t(x)$ is the $t$-fraction of the geodesic $\alpha_x$. Then it follows that $\gamma_0(x)=u_0(x)$, $\gamma_1(x)=u_1(x)$ and 
%\begin{align*}
%	d(\gamma_t(x),\gamma_s(x))=|t-s|d(u_0(x),u_1(x))
%\end{align*}
%for each $s,t\in [0,1]$ and for all $x\in X$. Consequently, we have
%\begin{align*}
%	D(\gamma_t,\gamma_s)=\Big(\int_Xd^2(\gamma_t(x),\gamma_s(x))d\mu(x) \Big)^{\frac{1}{2}}=|t-s|\Big(\int_Xd^2(u_0(x),u_1(x))d\mu(x) \Big)^{\frac{1}{2}},
%\end{align*}
%from which we deduce
%\begin{align*}
%	l_D(\gamma_t)=\Big(\int_Xd^2(u_0(x),u_1(x))d\mu(x) \Big)^{\frac{1}{2}}=D(u_0,u_1).
%\end{align*}
%It remains to show the triangle comparison. Fix $u_0,u_1,u_2\in L^2(X,Y)$ and let $u_t$ be the $t$-fraction of the geodesic connecting $u_1$ and $u_2$. Then for each $x\in X$ and $t\in [0,1]$, the NPC property of $Y$ gives
%\begin{align*}
% d^2(u_0(x),u_1(x))&\leq (1-t)d^2(u_0(x),u_1(x))+td^2(u_0(x),u_2(x))\\
% &\qquad -t(1-t)d^2(u_1(x),u_2(x)).
%\end{align*}
%Integrating with respect to $x\in X$, we conclude that
%\begin{align*}
%	D^2(u_0,u_t)\leq (1-t)D^2(u_0,u_t)+tD^2(u_0,u_2)-t(1-t)D^2(u_1(x),u_2(x)).
%\end{align*}
%This completes the proof.
%\end{proof}

\begin{proof}[Proof of Theorem~\ref{thm:gradient flow}]
	Since $u\mapsto E(u)$ is a lower semicontinuous convex functional on the NPC space $L^2(X,Y)$ (by Lemma~\ref{lemma:L2 is NPC}), the result follows immediately from~\cite[Theorem 1.13]{m98}.   
\end{proof}

\begin{proof}[Proof of Corollary~\ref{coro:gradient flow}]
	As the Sobolev energy of $u(t)$ decreases, it is uniformly bounded. As $X$ is compact, Theorem~\ref{thm:gradient flow} implies that $u(t)$ has bounded $L^2$-norm. By the precompactness of Sobolev spaces, there is subsequence $\{u(t_k)\}$ that converges in $L^2$ as $t_k\to \infty$. By~\cite[Proposition 2.4]{m98}, it will converges to an energy minimizer, which is a constant mapping.
\end{proof}

Fix an open set $\Omega\subset X$ and we set 
\begin{align*}
	\mathcal{KS}^{1,2}_\phi(\Omega,Y)=\Big\{u\in KS^{1,2}_\phi(\Omega,Y): E(u)\leq E(\phi) \Big\}.
\end{align*}

\bl\label{lemma:KS is NPC}
The metric space $\mathcal{KS}^{1,2}_\phi(\Omega,Y)$ (equipped with the metric $D$) is NPC.
\el 
\begin{proof}
	For each sequence $\{u_k\}\subset \mathcal{KS}^{1,2}_\phi(\Omega,Y)$ with $u_k\to u$ in $L^2(X,Y)$, we have $E(u)\leq E(\phi)$ by the lower semicontinuity of the Dirichlet energy. 
	%It remains to show that $u\in$Since $X$ is , we have by Theorem~\ref{thm:KS is included in Newton} that
	Note that
	$$E(u_k)+E(u)\leq 2E(\phi),$$ 
	which is uniformly bounded and $u_k\to u$ in $L^2(X,Y)$. We may apply Lemma~\ref{lemma:key lemma for existence} to conclude that $u\in \mathcal{KS}^{1,2}_\phi(\Omega,Y)$. It remains to show that $\mathcal{KS}^{1,2}_\phi(\Omega,Y)$ is a length space, which suffices as $\mathcal{KS}^{1,2}_\phi(\Omega,Y)$ then inherits the NPC property from $L^2(X,Y)$. However, this is clear since $\mathcal{KS}^{1,2}_\phi(\Omega,Y)$ is a convex subset of $L^2(X,Y)$ (due to the convexity of the Korevaar-Schoen energy).

\end{proof}

It is now possible to apply~\cite[Theorem 1.13]{m98}. Note that the inequality $E(u)\leq E(\phi)$ in the definition of $\mathcal{KS}^{1,2}_\phi(\Omega,Y)$ is only a technical requirement for showing $\mathcal{KS}^{1,2}_\phi(\Omega,Y)$ is closed. It has no effect on the flow as each time the energy decreases. When $\Omega$ is relatively compact in $X$, by the proof of Theorem~\ref{thm:main thm existence and uniqueness}, any minimizing sequence converges to the unique minimizer for the Dirichlet energy of Korevaar and Schoen. Thus we obtain the following theorem via~\cite[Theorem 1.13]{m98}.

\begin{proof}[Proof of Theorem~\ref{thm:initial boundary value problem}]
	This is a direct consequence of Lemma~\ref{lemma:KS is NPC} and~\cite[Theorem 1.13]{m98}.
\end{proof}

%\section*{Appendix}

\subsection*{Acknowledgements}

The author would like to thank Prof.~Yau Shing-Tung for his insightful comments on the theory of harmonic mappings and for his encouragement on developing the theory in the singular space setting during the 7th ICCM in Beijing 2016, which is the main motivation of this paper. He is grateful to Prof.~Jost for kindly sharing his papers~\cite{j94,j97} and for his valuable comments on the background of harmonic mappings, to Prof.~Chen for kindly sharing his paper~\cite{c95}, and to Profs.~Zhang and Zhu for their interest in this work and for kindly sharing the references~\cite{h16,hz16,zz16}. 
%We thank the anonymous referees for their careful reading and useful comments that greatly improve the presentation.

%\bibliography{bibli-Rellich}

\end{document}